\documentclass[preprint,12pt]{elsarticle}

\usepackage{amssymb}
\usepackage{amsthm}
\usepackage{amsmath}
\usepackage{fancybox}

\newtheorem{thrm}{Theorem}

\newtheorem{thm}{Theorem}
\newtheorem{lem}{Lemma}
\newdefinition{rem}{Remark}
\newdefinition{exa}{Example}
\newproof{pf}{{\bf Proof}}
\newproof{pfadd}{{\bf Proof (of additivity of~$\riV$)}}
\newproof{pflsc}{{\bf Proof (of lower semicontinuity of~$\riV$)}}
\newproof{plim}{{\bf Proof (of \protect\eq{e:ttos} for~$\riV$)}}
\newproof{pfBV}{{\bf Proof of Theorem~\protect\ref{t:BV}}}
\newproof{pfA}{{\bf Proof (of Theorem~\protect\ref{t:A})}}

\newproof{pfirst}{{\bf Proof of Theorem~\protect\ref{t:1}}}
\newproof{pscnd}{{\bf Proof of Theorem~\protect\ref{t:2}}}
\newcommand{\eq}[1]{\mbox{\rm(\ref{#1})}}

\catcode`\@=11  \@addtoreset{equation}{section}  \catcode`\@=12

\newcommand{\Rb}{\mathbb{R}}  \newcommand{\Nb}{\mathbb{N}}
\newcommand{\Oc}{\mathcal{O}}
\newcommand{\es}{\varnothing}
\newcommand{\vep}{\varepsilon}
\newcommand{\nV}{\mbox{\rm V}}
\newcommand{\riV}{{\scriptstyle\overrightarrow{\mbox{\rm V}}}}
\newcommand{\leV}{{\scriptstyle\overleftarrow{\mbox{\rm V}}}}
\newcommand{\ori}[1]{{\scriptscriptstyle\overrightarrow{\mbox{$#1$}}}}
\newcommand{\ole}[1]{{\scriptscriptstyle\overleftarrow{\mbox{$#1$}}}}
\newcommand{\ov}[1]{\overline{#1}}

\journal{arXiv}

\begin{document}

\begin{frontmatter}

\title{Asymmetric variations of multifunctions\\
       with application to functional inclusions}

\author[hsenn]{Vyacheslav V.~Chistyakov\corref{cor1}}
\ead{czeslaw@mail.ru, vchistyakov@hse.ru}



\cortext[cor1]{Corresponding author.}

\address[hsenn]{Department of Informatics, Mathematics and Computer Science,\\
National Research University Higher School of Economics,\\
Bol'shaya Pech{\"e}rskaya Street 25/12, Nizhny Novgorod\\
603155, Russia}

\begin{abstract}
Under certain initial conditions, we prove the existence of set-valued selectors of
univariate compact-valued multifunctions of bounded (Jordan) variation when the notion
of variation is defined taking into account only the Pompeiu asymmetric excess between
compact sets from the target metric space. For this, we study subtle properties of
the directional variations. We show by examples that all assumptions in the main
existence result are essential. As an application, we establish the existence of
set-valued solutions $X(t)$ of bounded variation to the functional inclusion of the
form $X(t)\subset F(t,X(t))$ satisfying the initial condition $X(t_0)=X_0$. 
\end{abstract}

\end{frontmatter}

\section{Introduction} \label{s:intro}

This paper is devoted to the existence of (set-valued) selectors (or selections) with
prescribed initial conditions of multifunctions (=set-valued mappings) with compact
images from a metric space. Continuous and Lipschitz continuous selectors exist,
in general, for multifunctions with convex (and closed) images from a Banach space
(\cite{AF,Michael,RS}). Many examples are known when a continuous (or even
H{\"o}lder continuous) multifunction from a closed interval $T=[a,b]\subset\Rb$
into a family of compact subsets of a ball in $\Rb^2$, or a Lipschitz continuous multifunction
from $\Rb^3$ into a family of compact subsets of a ball in $\Rb^3$ have no 
continuous selector
(\cite{AC,Posit98,Hermes69,Kupka96,Michael}).

In contrast to this, it was shown in \cite{Hermes} that a Lipschitz continuous
multifunction $F$ on $[a,b]$ with compact images from $\Rb^N$ admits a Lipschitz
continuous selector, whose Lipschitz constant does not exceed that of the multifunction;
furthermore, if $F$ is only continuous and of bounded (Jordan) variation, then
it admits a \emph{continuous\/} selector. Similar assertions for Lipschitz and
absolutely continuous multifunctions with convex and nonconvex images from $\Rb^N$
were established in \cite{GurKos,Kik-Tom,Kupka,Zhu}. The selector results
from \cite{Hermes} were extended in \cite[Suppl.~1]{Mor} for multifunctions
$F$ on $[a,b]$ with compact graphs and images in a Banach space, and in \cite{Slezak}
for metric space valued multifunctions with compact images. It is to be noted that the
compactness arguments in these references were based on Arzel{\`a}-Ascoli's theorem.

Changing the compactness arguments to generalized Helly's pointwise selection
principle, it was proved in \cite{JDCS} that a multifunction $F$ on $[a,b]$ of bounded
Jordan variation with compact images from a Banach space admits a selector,
whose total Jordan variation does not exceed that of $F$ and passes through a
given point in the graph of~$F$. Thus, ``nice'' selectors of compact-valued multifunctions
inherit boundedness of variation rather than continuity (\cite{Sovae}).
The results of \cite{JDCS} were then extended to multifunctions with compact
images in a metric space and having certain regularity such as absolute continuity,
Lipschitz continuity, boundedness of Riesz-Orlicz-Medvedev generalized variation,
boundedness of essential variation, measurability in the first variable and boundedness
of variation in the second variable (\cite{JMAA,MatSb,Pontr,Posit01,MZ,%
Sovae,JFA05,Studia}). The existence of selections of bounded variation was applied
to the study of set-valued measure differential problems~(\cite{Satco}).

The notion of bounded variation of a multifunction $F$ from a subset $T\subset\Rb$
into a family of compact subsets of a metric space relies on the linear order in 
the domain $T$ of $F$ and the Pompeiu-Hausdorff metric in the range of $F$,
which is the maximum of two asymmetric excesses (cf.~Section~\ref{s:mr}).
In this paper, we prove the existence of set-valued selectors (in particular,
single-valued selectors) of bounded Jordan variation (with respect to the
Pompeiu-Hausdorff metric) under milder assumptions of boundedness of directional
variations of $F$, to the right or to the left, with respect to the Pompeiu excess only
(\cite{JMAA07} and Section~\ref{s:mr}). Our main result, Theorem~\ref{t:BV} in
Section~\ref{s:mr}, extends the results of \cite{JMAA07} for single-valued
selectors, and some partial results of \cite{Tretya} concerning set-valued selectors.

This paper, along with \cite{JMAA07}, may be considered as part of analysis in 
asymmetric spaces, giving an intuition of the notion of bounded variation in 
the case when the distance function on the underlying space does not satisfy the 
symmetry axiom (for asymmetric analysis in normed spaces, cf.~\cite{Cobzas}).

The paper is organized as follows. In Section~\ref{s:mr}, we review certain definitions
and facts and present our main result, Theorem~\ref{t:BV}. Section~\ref{s:dirvar}
is devoted to the (subtle) properties of directional variations, presented for
compact-valued multifunctions defined on an arbitrary subset of~$\Rb$ (hence
suitable for linearly ordered sets as the domains). In Section~\ref{s:pro}, we prove 
Theorem~\ref{t:BV} and exhibit some of its consequences. We illustrate our result
by suitable examples in Section~\ref{s:exa}. In Section~\ref{s:Func_Em}, we prove
the existence of set-valued solutions of bounded variation to the functional inclusion
of the form $X(t)\!\subset\! F(t,X(t))$~for~all~$t\!\in\! T$.

\section{Main Result} \label{s:mr}

We begin by reviewing certain definitions and facts needed for our results.

Throughout the paper $(M,d)$ is a metric space with metric $d$.

Given two nonempty sets $X,Y\subset M$,
the \emph{Pompeiu excess\/} ({\it {\'e}cart}, in French) {\it of $X$ over $Y$}
 is defined by (\cite[pp.~281--282]{Pomp}, \cite[Chapter~II]{CV})
  \begin{equation*}
e(X,Y)\equiv e_d(X,Y):=\sup_{x\in X}d(x,Y)=\inf\{r>0:X\subset\Oc_r(Y)\},
  \end{equation*}
where $d(x,Y):=\inf_{y\in Y}d(x,y)$ is the distance\footnote{$d(x,\es):=+\infty.$}
from $x\in M$ to $Y$,
and the set $\Oc_r(Y):=\{x\in M:d(x,Y)<r\}$ is the open $r$-neighbourhood%
\footnote{$\Oc_r(\es)\!:=\!\es$ for $r>0$; $e(\es,Y):=0$ for any $Y\subset M$,
and $e(X,\es):=+\infty$ if $\es\!\ne\! X\!\subset\! M$.}
of~$Y$, $r>0$ (if $Y=\{y\}$,
$\Oc_r(y)\equiv\Oc_r(Y)$ is simply the open ball of radius $r$ centered at $y\in M$).
Note that $e(X,Y)\ne e(Y,X)$ in general.

The well-known properties of $e(\cdot,\cdot)$ are as follows. Given $X,Y,Z\subset M$:
\par\vspace{2pt}
(a) $e(X,Y)\!=\!0$ iff\footnote{`iff' means as usual `if and only if'.} $X\!\subset\!\ov{Y}$,
where $\ov{Y}\!=\!\bigcap_{r>0}\Oc_r(Y)$ is the closure of $Y$~in~$M$;
\par\vspace{2pt}
(b) $e(X,Y)\le e(X,Z)+e(Z,Y)$ (triangle inequality for~$e$);
\par\vspace{2pt}
(c) $e(X,Y)<+\infty$ provided $X$ and $Y$ are bounded (and, in particular bounded
and closed, or compact).

The \emph{Pompeiu-Hausdorff distance\/} between two sets $X,Y\subset M$
 is defined by (e.g., \cite{Pomp}, \cite[Chapter~II]{CV})%
\footnote{Pompeiu \cite{Pomp} symmetrized the excess $e$ by
 $d_P(X,Y)\!:=\!e(X,Y)\!+\!e(Y,X)$ (cf. also~\cite{Ber-Pac}).}
  \begin{equation*}
d_H(X,Y)\!:=\max\{e(X,Y),e(Y,X)\}\!=\!
\inf\{r\!>\!0:\mbox{$X\!\subset\!\Oc_r(Y)$\,and\,$Y\!\subset\!\Oc_r(X)$}\}.
  \end{equation*}
It follows from (a), (b), and (c) above that $d_H$ is a \emph{metric\/}
(with finite values), called the \emph{Pompeiu-Hausdorff metric}, on the family
$\mbox{bcl}(M)$ of all nonempty \emph{bounded closed\/} subsets of $M$ and, in 
particular on the family $\mbox{c}(M)$ of all nonempty \emph{compact\/} subsets of~$M$.

A \emph{multifunction\/} (or \emph{set-valued mapping\/}) from a nonempty set $T$
into $M$ is a rule $F$, which assigns to each $t\in T$ a unique subset $F(t)\subset M$;
in symbols, $F:T\to\mathcal{P}(M)$, where $\mathcal{P}(M)$ is the power set of $M$
(=\,the family of all subsets of~$M$). A multifunction $\Gamma:T\to\mathcal{P}(M)$
is said to be a \emph{set-valued selector\/} of $F:T\to\mathcal{P}(M)$ on $T$ provided
$\Gamma(t)\subset F(t)$ for all $t\in T$. Clearly, $F$ is a set-valued selector of itself.
If $\Gamma:T\to M$ is single-valued and $\Gamma(t)\in F(t)$ for all $t\in T$, then
$\Gamma$ is called a \emph{selector\/} (or \emph{selection\/}) of~$F$ on~$T$.

Of main interest in this paper are multifunctions $F:T\to\mbox{c}(M)$
with $T$ a nonempty subset of the reals $\Rb$. Such an $F$ is said to be
\emph{of bounded variation\/} (with respect to~$d_H$) provided its (total)
 \emph{Jordan variation\/}
  \begin{equation*}
\nV(F,T)\!:=\sup_\pi\sum_{i=1}^m d_H\bigl(F(t_{i-1}),F(t_i)\bigr)\quad
\mbox{is finite\quad ($\nV(F,\es):=0$),}
  \end{equation*}
the supremum being taken over all partitions $\pi$ of the set $T\subset\Rb$, i.e.,
$m\in\Nb$ and $\pi=\{t_i\}_{i=0}^m\subset T$ such that $t_{i-1}\le t_i$
for all $i\in\{1,\dots,m\}$.

The following theorem on the existence of set-valued selectors of bounded variation
was established in \cite[Theorems~10.1 and~5.1]{Sovae}:

\begin{thrm} \label{t:A}
Given $T\subset\Rb$, $t_0\in T$, $X_0\in\mbox{\rm c}(M)$, and
$F:T\to\mbox{\rm c}(M)$ such that $\nV(F,T)<+\infty$, there exists
a set-valued selector $\Gamma:T\to\mbox{\rm c}(M)$ of $F$ on $T$ such that
$d_H(X_0,\Gamma(t_0))\le e(X_0,F(t_0))$ and\/ $\nV(\Gamma,T)\le\nV(F,T)$.
{\rm(}In particular, if $X_0\subset F(t_0)$, the first inequality above gives
$\Gamma(t_0)=X_0$.{\rm)} Furthermore, if $x_0\in M$ and $X_0=\{x_0\}$,
then $\Gamma:T\to M$ may be chosen to be single-valued
and such that $d(x_0,\Gamma(t_0))\le d(x_0,F(t_0))$.
\end{thrm}

Theorem~\ref{t:A} already contains as particular cases many previously known results
\cite{JMAA}, \cite{JDCS}--\cite{MZ}, \cite{Studia,GurKos,Hermes,Mor,Slezak},
concerning single-valued selectors.

The purpose of this paper is to drop the assumption $\nV(F,T)\!<\!+\infty$, replacing it by
`directional' ones $\riV(F,T)\!<\!+\infty$ or $\leV(F,T)\!<\!+\infty$ (see below), which,
as will be shown, still guarantees the existence of set-valued selectors $\Gamma$ of $F$
on $T$ of bounded variation (with respect to~$d_H$).

In order to do it, given $F:T\to\mbox{c}(M)$, following \cite{JMAA07} the quantities
  \begin{equation*}
\riV(F,T)\!:=\!\sup_\pi\sum_{i=1}^m e\bigl(F(t_{i-1}\!),F(t_i)\bigr)
\mbox{ and \,}
\leV(F,T)\!:=\!\sup_\pi\sum_{i=1}^m e\bigl(F(t_i),F(t_{i-1}\!)\bigr)
  \end{equation*}
are said to be the \emph{directional\/} (or \emph{excess\/}) \emph{variations of $F$
to the right\/} and \emph{to the left}, respectively ($\riV(F,\es)=\leV(F,\es):=0$). Clearly,
  \begin{equation*}
\max\bigl\{\riV(F,T),\leV(F,T)\bigr\}\le\nV(F,T)\le\riV(F,T)+\leV(F,T),
  \end{equation*}
and so, $\nV(F,T)$ is finite iff both $\riV(F,T)$ and $\leV(F,T)$ are finite.
Note that if $F:T\to M$ is single-valued, the quantity $\nV(F,T)=\riV(F,T)=\leV(F,T)$ 
is the usual (Jordan) variation of $F$ on $T$ (e.g., \cite{MatSb}).

Our first main result, extending Theorem~\ref{t:A}, is the following theorem
on the existence of set-valued selectors of bounded variation.
For brevity, we write $T_S:=T\cap S$ for $S\subset\Rb$ (e.g.,
$T_{[t_0,+\infty)}=T\cap[t_0,+\infty)$ for $t_0\in T$, etc.).

\begin{thm} \label{t:BV}
Suppose $T\subset\Rb$, $t_0\in T$, $X_0\in\mbox{\rm c}(M)$, and
$F:T\to\mbox{\rm c}(M)$.

  \begin{itemize}
\item[{\rm(a)}] Let $a:=\inf T\in T$. If\/ $\riV(F,T)<+\infty$, then there exists
a set-valued selector of bounded variation $\Gamma:T\to\mbox{\rm c}(M)$
of $F$ on $T$ such that  $d_H(X_0,\Gamma(t_0))\le e(X_0,F(t_0))$,
$\nV(\Gamma,T_{[a,t_0)})\le\riV(F,T_{[a,t_0)})$,
$\nV(\Gamma,T_{[t_0,+\infty)})$ $\le\riV(F,T_{[t_0,+\infty)})$,
and
  \begin{equation} \label{e:s1t0}
\nV(\Gamma,T)-J_a(\Gamma,t_0)\le
  \riV(F,T_{[a,t_0)})+\riV(F,T_{[t_0,+\infty)})\le\riV(F,T),
  \end{equation}
where $J_a(\Gamma,t_0):=0$ if $t_0=a$, and if $t_0>a$ and $s_0:=\sup T_{[a,t_0)}$,
  \begin{equation} \label{e:J}
J_a(\Gamma,t_0):=\left\{
  \begin{array}{ccc}
d_H(\Gamma(s_0),\Gamma(t_0)), & \,\,\mbox{if}\,\,\, & s_0\in T_{[a,t_0)},\\[4pt]
\displaystyle\lim_{T\ni t\to s_0-0}\!\nV(\Gamma,T_{[t,t_0]}),
  & \,\,\mbox{if}\,\,\, & s_0\notin T_{[a,t_0)}.
  \end{array}\right.
  \end{equation}
In particular, if $s_0=t_0$, then $J_a(\Gamma,t_0)=\displaystyle%
\lim_{T\ni t\to t_0-0}\!d_H(\Gamma(t),\Gamma(t_0))$.

\item[{\rm(b)}] Let $b:=\sup T\in T$. If\/ $\leV(F,T)<+\infty$, then there exists
a set-valued selector of bounded variation $\Gamma:T\to\mbox{\rm c}(M)$
of $F$ on $T$ such that  $d_H(X_0,\Gamma(t_0))\le e(X_0,F(t_0))$,
$\nV(\Gamma,T_{(t_0,b]})\le\leV(F,T_{(t_0,b]})$,
$\nV(\Gamma,T_{(-\infty,t_0]})$ $\le\leV(F,T_{(-\infty,t_0]})$,
and
  \begin{equation*}
\nV(\Gamma,T)-J_b(\Gamma,t_0)\le
  \leV(F,T_{(t_0,b]})+\leV(F,T_{(-\infty,t_0]})\le\leV(F,T),
  \end{equation*}
where $J_b(\Gamma,t_0):=0$ if $t_0=b$, and if $t_0<b$ and $s_0:=\inf T_{(t_0,b]}$,
  \begin{equation*}
J_b(\Gamma,t_0):=\left\{
  \begin{array}{ccc}
d_H(\Gamma(t_0),\Gamma(s_0)), & \,\,\mbox{if}\,\,\, & s_0\in T_{(t_0,b]},\\[4pt]
\displaystyle\lim_{T\ni t\to s_0+0}\!\nV(\Gamma,T_{[t_0,t]}),
  & \,\,\mbox{if}\,\,\, & s_0\notin T_{(t_0,b]}.
  \end{array}\right.
  \end{equation*}
In particular, if $s_0=t_0$, then $J_b(\Gamma,t_0)=\displaystyle%
\lim_{T\ni t\to t_0+0}\!d_H(\Gamma(t_0),\Gamma(t))$.
  \end{itemize}
\end{thm}

This theorem will be proved in Section~\ref{s:pro}, where it will also be shown that
Theorem~\ref{t:BV} implies Theorem~\ref{t:A}.
A somewhat free and intuitive interpretation of Theorem~\ref{t:BV} can be given as
follows. We may imagine $F$ to be a road with the value $F(t)$ as a section at a given
coordinate~$t$. If $t$ increases, the section $F(t)$ moves in one direction, say, to the
right. Theorem~\ref{t:BV} asserts that if the road $F$ is properly built/controlled
(i.e., its variation to the right is finite), then an extended object (e.g., a car) can
freely pass it.

\section{Properties of Directional Variations} \label{s:dirvar}

In this section, we gather auxiliary facts needed for the proof of Theorem~\ref{t:BV}.
Let $T\subset\Rb$ and $F:T\to\mbox{c}(M)$ be a multifunction.

\subsection{Monotonicity of $F$} \label{ss:mon}

By the definition of $\nV$, $\nV(F,T)=0$ iff $F$ is constant on $T$. The definition of~$\riV$
and property (a) of $e$ in Section~\ref{s:mr} imply $\riV(F,T)=0$ iff $F$ is
\emph{nondecreasing\/} on $T$ in the sense that $F(s)\subset F(t)$ for all $s,t\in T$
with $s\le t$. Similarly, $\leV(F,T)=0$ iff $F$ is \emph{nonincresing\/} on $T$, i.e.,
$F(s)\supset F(t)$ for all $s,t\in T$,~$s\le t$.

\subsection{Additivity of Variations}

Since $(\mbox{c}(M),d_H)$ is a metric space, it is known (e.g., \cite[2.19]{Var},
\cite{JDCS}, \cite{MatSb}) that $\nV$ is \emph{additive\/} (in the second variable)
in the sense that
  \begin{equation} \label{e:Vadd}
\nV(F,T)=\nV(F,T_{(-\infty,t]})+\nV(F,T_{[t,+\infty)})\quad\mbox{for all \,\,$t\in T$.}
  \end{equation}
We assert that the additivity property \eq{e:Vadd} holds also for $\riV$ and $\leV$
in place of~$\nV$. Since the Hausdorff excess $e$ is not symmetric, we have to take
care of the order of its arguments. So, we explicitly verify \eq{e:Vadd} at least for~$\riV$.

\begin{pfadd}
Given $\xi,\eta\in\Rb$ such that $\xi<\riV(F,T_{(-\infty,t]})$ and
$\eta<\riV(F,T_{[t,+\infty)})$, there are $m,n\in\Nb$, a partition $\{t_i\}_{i=0}^m$
of $T_{(-\infty,t]}$ with $t_m=t$ and a partition $\{s_j\}_{j=0}^n$ of $T_{[t,+\infty)}$
with $s_0=t$ such that
  \begin{equation} \label{e:xi1}
\xi<\sum_{i=1}^m e(F(t_{i-1}),F(t_i))\quad\mbox{and}\quad
\eta<\sum_{j=1}^n e(F(s_{j-1}),F(s_j)).
  \end{equation}
Since $\{t_i\}_{i=0}^m\cup\{s_j\}_{j=0}^n$ is a partition of $T$, we get, from \eq{e:xi1},
$\xi+\eta<\riV(F,T)$, which, due to the arbitrariness of $\xi$ and $\eta$ as above,
implies the inequality
  \begin{equation*}
\riV(F,T_{(-\infty,t]})+\riV(F,T_{[t,+\infty)})\le\riV(F,T).
  \end{equation*}

In order to prove the reverse inequality, let $\xi\in\Rb$ be arbitrary such that
$\xi<\riV(F,T)$. Then, there are $m\in\Nb$ and a partition $\pi=\{t_i\}_{i=0}^m$
of $T$ such that the first inequality in \eq{e:xi1} holds. If $t_m\le t$, then $\pi$ is a
partition of $T_{(-\infty,t]}$, and if $t\le t_0$, then $\pi$ is a partition of $T_{[t,+\infty)}$,
and so, in the either case, the first inequality in \eq{e:xi1} implies
  \begin{equation} \label{e:x2}
\xi<\riV(F,T_{(-\infty,t]})+\riV(F,T_{[t,+\infty)}).
  \end{equation}
Now, suppose $t_{k-1}\le t\le t_k$ for some $k\in\{1,\dots,m\}$. By the triangle
inequality for~$e$, we have
  \begin{equation} \label{e:kk}
e(F(t_{k-1}),F(t_k))\le e(F(t_{k-1}),F(t))+e(F(t),F(t_k)).
  \end{equation}
Since $\{t_i\}_{i=0}^{k-1}\cup\{t\}$ is a partition of $T_{(-\infty,t]}$ and
$\{t\}\cup\{t_i\}_{i=k}^m$ is a partition of $T_{[t,+\infty)}$, from the first inequality
in \eq{e:xi1} and \eq{e:kk}, once again we get \eq{e:x2}. By the arbitrariness
 of $\xi<\riV(F,T)$,
\eq{e:x2} implies
   \begin{equation*}
\riV(F,T)\le\riV(F,T_{(-\infty,t]})+\riV(F,T_{[t,+\infty)}),
  \end{equation*}
 which was to be proved.
\qed\end{pfadd}

\subsection{Bounded Directional Variations} \label{ss:bdv}

Recall that the function $v_F(t)\!:=\!\nV(F,T_{(-\infty,t]})$, $t\!\in\! T$, is said to be
the {\it var\-i\-a\-tion function\/} of $F$ on~$T$. 
We define the \emph{variation function\/} of $F$ to the right (to the left) by
$\ori{v}_F(t):=\riV(F,T_{(-\infty,t]})$ (by $\ole{v}_F(t):=\leV(F,T_{(-\infty,t]})$,
respectively) for all $t\in T$. Clearly, $v_F,\ori{v}_F,\ole{v}_F:T\to[0,+\infty]$ are
\emph{nondecreasing\/} on~$T$. If $F$ is clear from the context, we omit the
subscript $F$ and write $v$, $\ori{v}$, and $\ole{v}$, respectively.

The following characterization holds for multifunctions of bounded directional
variation (which will be useful in Section~\ref{s:Func_Em}).

\begin{lem} \label{l:char} Given $F:T\to\mbox{\rm c}(M)$, $\riV(F,T)<+\infty$ iff
there is a nondecreasing bounded function $\varphi:T\to\Rb$ such that
$e(F(s),F(t))\le\varphi(t)-\varphi(s)$ for all $s,t\in T$ with $s\le t;$ moreover,
$\riV(F,T)\le\nV(\varphi,T)$. {\rm(}A similar assertion holds for $\leV(F,T)<+\infty$ if
the inequality is written as $e(F(t),F(s))\le\varphi(t)-\varphi(s)$.{\rm)}
\end{lem}

\begin{pf}
($\Rightarrow\!$) Let $\ori{v}(t)=\riV(F,T_{(-\infty,t]})$, $t\in T$, be the variation
function (to the right) of~$F$. It is nondecreasing on $T$ and bounded:
$\sup_{t\in T}\ori{v}(t)\le\riV(F,T)$. Given $s,t\in T$, $s\le t$, the additivity
of $\riV$ implies
  \begin{equation*}
e(F(s),F(t))\le\riV(F,T_{[s,t]})=\riV(F,T_{(-\infty,t]})-\riV(F,T_{(-\infty,s]})
=\ori{v}(t)-\ori{v}(s).
  \end{equation*}
It remains to set $\varphi:=\ori{v}$.

($\!\Leftarrow$) Let $m\in\Nb$ and $\pi=\{t_i\}_{i=0}^m$ be a partition of $T$.
Since $t_{i-1}\le t_i$ for all $i=1,\dots,m$, by the assumption we get
  \begin{align*}
\sum_{i=1}^me(F(t_{i-1}),F(t_i))&\le\sum_{i=1}^m\bigl(\varphi(t_i)-\varphi(t_{i-1})\bigr)
  =\varphi(t_m)-\varphi(t_0)\\
&\le\sup_{t\in T}\varphi(t)-\inf_{s\in T}\varphi(s)=\nV(\varphi,T)<+\infty.
  \end{align*}
It remains to take the supremum over all partitions $\pi$ of $T$.
\qed\end{pf}

\begin{rem} \label{p:7}
It is known (\cite[1.23]{Var}, \cite[Lemma~11]{MZ}) that if $F:T\to\mbox{c}(M)$
and $\nV(F,T)\!<\!+\infty$, then the image $F(T):=\bigcup_{t\in T}F(t)$ is a
\emph{totally bounded\/} (hence separable) subset of~$M$ (if, in addition, $(M,d)$ is
complete, then the closure of $F(T)$ in $M$ is compact). However,
condition $\nV(F,T)\!<\!+\infty$ cannot be replaced neither by
$\riV(F,T)\!<\!+\infty$ nor by $\leV(F,T)\!<\!+\infty$ (see Example~\ref{ex:5.5}).
\end{rem}

\subsection{Lower Semicontinuity of Variations} \label{ss:lsc}

If a sequence $\{F_n\}_{n=1}^\infty$ of multifunctions $F_n:T\to\mbox{c}(M)$
\emph{converges in $\mbox{\rm c}(M)$ pointwise on $T$} to $F$ (i.e.,
$d_H(F_n(t),F(t))\to0$ as $n\to\infty$ for all $t\in T$), then
  \begin{equation} \label{e:Vsem}
\nV(F,T)\le\liminf_{n\to\infty}\nV(F_n,T).
  \end{equation}
This property is known as the (sequential) \emph{lower semicontinuity\/} of~$\nV$
(in the first variable) for metric space valued functions; cf. \cite[Proposition~2.1(V7)]{MatSb}.

\smallbreak
We assert that \eq{e:Vsem} is valid for $\riV$ and $\leV$ in place of $\nV$, as well.

\begin{pflsc}
For any $n\in\Nb$ and $s,t\in T$, the triangle inequality for $e$ implies
  \begin{equation*}
|e(F_n(s),F_n(t))-e(F(s),F(t))|\le d_H(F_n(s),F(s))+d_H(F_n(t),F(t)),
  \end{equation*}
and so, by the pointwise convergence of $F_n$ to $F$, $e(F_n(s),F_n(t))\!\to\! e(F(s),F(t))$
as $n\to\infty$. Given $m\in\Nb$ and a partition $\pi=\{t_i\}_{i=0}^m$ of $T$,
by definition of~$\riV$,
  \begin{equation*}
\sum_{i=1}^m e(F_n(t_{i-1}),F_n(t_i))\le\riV(F_n,T)\quad\mbox{for all \,\,$n\in\Nb$.}
  \end{equation*}
Passing to the limit inferior as $n\to\infty$, we get
  \begin{equation*}
\sum_{i=1}^m e(F(t_{i-1}),F(t_i))\le\liminf_{n\to\infty}\riV(F_n,T),
  \end{equation*}
and it remains to take the supremum over all partitions $\pi$ of $T$.
\qed\end{pflsc}

\subsection{Limit Equalities for the Variations}

The following equalities are known (\cite{Var,JDCS}) for (multi)functions $F$ on $T$
with values in a metric space (in particular in $\mbox{c}(M)$):
  \begin{align}
&\mbox{if $s\!=\!\sup T\!\in\!\Rb\cup\{+\infty\}$ and $s\!\notin\! T$, then
  $\displaystyle\nV(F,T)\!=\!\lim_{T\ni t\to s}\nV(F,T_{(-\infty,t]});$} \label{e:ttos}\\
&\mbox{if \,$i\!=\!\inf T\!\in\Rb\cup\{-\infty\}$ \,and $i\!\notin\! T$, then
  $\displaystyle\nV(F,T)\!=\!\lim_{T\ni t\to i}\nV(F,T_{[t,+\infty)}).$}\label{e:ttoi}
  \end{align}

We are going to show that these assertions hold for $\riV$ and $\leV$ as well.

\begin{plim}
Since $\ori{v}$ is nondecreasing on $T$, the limit on the right in \eq{e:ttos} (with $\nV$
replaced by~$\riV$) exists in $[0,+\infty]$ and is equal to $\sup_{t\in T}\ori{v}(t)$.
Given $t\in T$, we have $T_{(-\infty,t]}\subset T$, and so, $\ori{v}(t)\le\riV(F,T)$,
which implies $\lim_{T\ni t\to s}\ori{v}(t)\le\riV(F,T)$. Conversely, given $\xi\in\Rb$
with $\xi<\riV(F,T)$, there are $m\in\Nb$ and a partition $\pi=\{t_i\}_{i=0}^m$ of $T$
(and so, \mbox{$t_0\le t_1\le\dots\le t_{m-1}\le t_m<s$}) such that
  \begin{equation*}
\xi<\sum_{i=1}^m e(F(t_{i-1}),F(t_i))\le\ori{v}(t_m)\le\lim_{T\ni t\to s}\ori{v}(t).
  \end{equation*}
Passing to the limit as $\xi\to\riV(F,T)$, we get $\riV(F,T)\le\lim_{T\ni t\to s}\ori{v}(t)$.
\qed\end{plim}

\subsection{Jump Formulas for the Variations}

We say that $t\in\Rb$ is a \emph{left\/} (\emph{right\/}) \emph{limit point\/} of~$T$ if
$T_{(t-\vep,t)}\ne\es$ (respectively, $T_{(t,t+\vep)}\ne\es$) for all $\vep>0$.
A point $t\in T$, which is not a left (or right) limit point of $T$, is called the \emph{left\/}
(or \emph{right\/}) \emph{isolated point\/} of~$T$.

By virtue of \cite[Lemma~4.2]{MatSb}, the following two assertions hold.

\begin{itemize}
\item[{\rm(a)}] If $t\in T$ is a left limit point of $T$ and $\nV(F,T)<+\infty$, then
  \begin{equation} \label{e:VFtm}
v(t)=v(t-0)+\lim_{T\ni s\to t-0}d_H(F(s),F(t)).%
\footnote{$v(t-0):=\lim_{T\ni s\to t-0}v(s)$ is the left limit of $v$ at $t$ along $T$,
i.e., $v(t-0)$ is the limit of $v(s)$ as $s\to t$ with $s\in T_{(-\infty,t)}$.}
  \end{equation}
\item[{\rm(b)}] If $t\in T$ is a right limit point of $T$ and $\nV(F,T)<+\infty$, then
  \begin{equation} \label{e:VF+0}
v(t+0)=v(t)+\lim_{T\ni s\to t+0}d_H(F(s),F(t)).%
\footnote{$v(t+0):=\lim_{T\ni s\to t+0}v(s)=%
  \lim_{s\to t,s\in T_{(t,+\infty)}}v(s)$ is the right limit of $F$ at~$t$.}
  \end{equation}
\end{itemize}

The second limit in the right-hand side of \eq{e:VFtm} (of \eq{e:VF+0}) is known as
the \emph{left\/} (\emph{right\/}, respectively) \emph{jump\/} of $F$ at~$t$.

Since $d_H$ is symmetric, the order of arguments $F(s)$ and $F(t)$ in \eq{e:VFtm}
and \eq{e:VF+0} does not matter. In the counterparts of (a) and (b) for $\riV$ and
$\leV$ below, we take care of the order of arguments in the excess $e(\cdot,\cdot)$.

\begin{lem} \label{l:jumps}
If $t\in T$ is a left limit point of $T$, then%
\footnote{The limit superior in ($\ori{\mbox{\rm a}}$) is the limit of
$\sup\{e(F(s),F(t)):s\in T_{(t-\vep,t)}\}$ as $\vep\to+0$;
in ($\ole{\mbox{\rm a}}$) the limit superior is understood similarly.}
  \begin{itemize}
\item[${\rm(\ori{\mbox{\rm a}})}$] $\ori{v}(t)=\ori{v}(t-0)+\displaystyle
  \limsup_{T\ni s\to t-0}e(F(s),F(t))$ \,if \,\,$\riV(F,T)<+\infty;$
\item[${\rm(\ole{\mbox{\rm a}})}$] $\ole{v}(t)=\ole{v}(t-0)+\displaystyle
  \limsup_{T\ni s\to t-0}e(F(t),F(s))$ \,if \,\,$\leV(F,T)<+\infty.$
  \end{itemize}

If $t\in T$ is a right limit point of $T$, then%
\footnote{The limit superior in ($\ori{\mbox{\rm b}}$) is the limit of
$\sup\{e(F(t),F(s)):s\in T_{(t,t+\vep)}\}$ as $\vep\to+0$;
the limit superior in ($\ole{\mbox{\rm b}}$) has a similar meaning.}
  \begin{itemize}
\item[${\rm(\ori{\mbox{\rm b}})}$] $\ori{v}(t+0)=\ori{v}(t)+\displaystyle
  \limsup_{T\ni s\to t+0}e(F(t),F(s))$ \,if \,\,$\riV(F,T)<+\infty;$
\item[${\rm(\ole{\mbox{\rm b}})}$] $\ole{v}(t+0)=\ole{v}(t)+\displaystyle
  \limsup_{T\ni s\to t+0}e(F(s),F(t))$ \,if \,\,$\leV(F,T)<+\infty.$
  \end{itemize}
\end{lem}

\begin{pf}
We concentrate on ($\ori{\mbox{\rm a}}$) and ($\ori{\mbox{\rm b}}$),
making only remarks on ($\ole{\mbox{\rm a}}$) and ($\ole{\mbox{\rm b}}$).

($\ori{\mbox{\rm a}}$)
First, we prove inequality $\ge$ in ($\ori{\mbox{\rm a}}$).
By the additivity of $\riV$, given $s\in T$ with $s<t$, we have
  \begin{equation*}
\ori{v}(s)+e(F(s),F(t))\le\ori{v}(s)+\riV(F,T_{[s,t]})=\ori{v}(t).
  \end{equation*}
Since $\ori{v}$ is nondecreasing and bounded on $T$ (by $\riV(F,T)$), and $t$ is a left
limit point of $T$, the left limit $\ori{v}(t-0)$ exists in $[0,+\infty)$ and is equal to
$\sup\{\ori{v}(s):s\in T,s<t\}$. Inequality $\ge$ in ($\ori{\mbox{\rm a}}$) follows
now from the properties of the limit superior as $T\!\ni\! s\!\to\! t\!-\!0$, since
the right-hand side in ($\ori{\mbox{\rm a}}$) is equal to
  \begin{equation*}
\lim_{T\ni s\to t-0}\!\ori{v}(s)\!+\!\limsup_{T\ni s\to t-0}\!e(F(s),F(t))\!=\!
\limsup_{T\ni s\to t-0}\!\bigl(\ori{v}(s)\!+\!e(F(s),F(t))\bigr)\le\ori{v}(t).
  \end{equation*}

Now we show that inequality $\le$ holds in ($\ori{\mbox{\rm a}}$). By the definition
of $\ori{v}(t)$, which is finite, for any $\xi\in\Rb$, $\xi<\ori{v}(t)$, there are $m\in\Nb$
and a partition $\pi=\{t_i\}_{i=0}^m\cup\{t\}$ of $T_{(-\infty,t]}$ with $t_m<t$
such that
  \begin{equation*}
\xi<\sum_{i=1}^m e(F(t_{i-1}),F(t_i))+e(F(t_m),F(t)).
  \end{equation*}
Given $s\in T_{(t_m,t)}$, the triangle inequality for $e$ implies
  \begin{equation*}
e(F(t_m),F(t))\le e(F(t_m),F(s))+e(F(s),F(t)),
  \end{equation*}
and so, $\xi<\ori{v}(s)+e(F(s),F(t))$. Passing to the limit superior as
$T\ni s\to t-0$, we get
  \begin{equation*}
\xi\le\ori{v}(t-0)+\limsup_{T\ni s\to t-0}e(F(s),F(t)).
  \end{equation*}
It remains to take into account the arbitrariness of $\xi<\ori{v}(t)$.

\smallbreak
($\ori{\mbox{\rm b}}$)
In order to prove inequality $\ge$ in ($\ori{\mbox{\rm b}}$), we make use of the
additivity of $\riV$ to get
  \begin{equation*}
\ori{v}(t)+e(F(t),F(s))\le\ori{v}(t)+\riV(F,T_{[t,s]})=\ori{v}(s)\quad
\mbox{for all $s\in T_{(t,+\infty)}$.}
  \end{equation*}
Now, it suffices to pass to the limit superior as $T\ni  s\to t+0$.

To show inequality $\le$ in ($\ori{\mbox{\rm b}}$), we apply the additivity property
of $\riV$ several times. To begin with, note that
  \begin{equation} \label{e:tin}
\riV(F,T_{[t,+\infty)})=\riV(F,T)-\ori{v}(t).
  \end{equation}
Given $\xi\in\Rb$ with $\xi<\riV(F,T_{[t,+\infty)})$, there exist $m\in\Nb$ and a
partition $\{t\}\cup\{t_i\}_{i=0}^m$ of $T_{[t,+\infty)}$ with $t<t_0$ such that
  \begin{equation*}
\xi<e(F(t),F(t_0))+\sum_{i=1}^m e(F(t_{i-1}),F(t_i)).
  \end{equation*}
Since, for any $s\!\in\! T_{(t,t_0)}$, $e(F(t),F(t_0))\!\le\! e(F(t),F(s))\!+\!e(F(s),F(t_0))$,
we~find
  \begin{align*}
\xi+\ori{v}(t_m)-\riV(F,T)\le\xi&<e(F(t),F(s))+\riV(F,T_{[s,t_m]})\\[2pt]
&=e(F(t),F(s))+\ori{v}(t_m)-\ori{v}(s),
  \end{align*}
and so, $\ori{v}(s)\le(\riV(F,T)-\xi)+e(F(t),F(s))$. Passing to the limit superior as
$T\ni s\to t+0$, we get
  \begin{equation*}
\ori{v}(t+0)\le\riV(F,T)-\xi+\limsup_{T\ni s\to t+0}e(F(t),F(s)).
  \end{equation*}
It remains to let $\xi$ tend to the value \eq{e:tin}.

\smallbreak
($\ole{\mbox{\rm a}}$), ($\ole{\mbox{\rm b}}$)
Here we note only that, by the additivity of $\leV$,
  \begin{align*}
\ole{v}(s)+e(F(t),F(s))\le\ole{v}(s)+\leV(F,T_{[s,t]})=\ole{v}(t)\quad
  \forall\,s\in T_{(-\infty,t)},\\[2pt]
\ole{v}(t)+e(F(s),F(t))\le\ole{v}(t)+\leV(F,T_{[t,s]})=\ole{v}(s)\quad
  \forall\,s\in T_{(t,+\infty)},
  \end{align*}
respectively.
\qed\end{pf}

As a corollary of \eq{e:ttos}, \eq{e:ttoi}, and Lemma~\ref{l:jumps}, we get the
following lemma, which is a generalization of Theorem~4.6 from~\cite{MatSb}.

\begin{lem} \label{l:opej}
If $t\in T$ is a left limit point of $T$, then
  \begin{itemize}
\item[${\rm(\ori{\mbox{\rm a}})}$]
$\riV(F,T_{(-\infty,t]})\!=\!\riV(F,T_{(-\infty,t)})\!+\!\displaystyle
  \limsup_{T\ni s\to t-0}e(F(s),F(t))$ \,if \,\,$\riV(F,T)\!<\!+\infty;$
\item[${\rm(\ole{\mbox{\rm a}})}$]
$\leV(F,T_{(-\infty,t]})\!=\!\leV(F,T_{(-\infty,t)})\!+\!\displaystyle
  \limsup_{T\ni s\to t-0}e(F(t),F(s))$ \,if \,\,$\leV(F,T)<+\infty.$
  \end{itemize}

If $t\in T$ is a right limit point of $T$, then
  \begin{itemize}
\item[${\rm(\ori{\mbox{\rm b}})}$]
$\riV(F,T_{[t,+\infty)})\!=\!\riV(F,T_{(t,+\infty)})\!+\!\displaystyle
  \limsup_{T\ni s\to t+0}e(F(t),F(s))$ \,if \,\,$\riV(F,T)\!<\!+\infty;$
\item[${\rm(\ole{\mbox{\rm b}})}$]
$\leV(F,T_{[t,+\infty)})\!=\!\leV(F,T_{(t,+\infty)})\!+\!\displaystyle
  \limsup_{T\ni s\to t+0}e(F(s),F(t))$ \,if \,\,$\leV(F,T)<+\infty.$
  \end{itemize}
\end{lem}

\begin{pf}
As in the proof of Lemma~\ref{l:jumps}, we concentrate 
on ($\ori{\mbox{\rm a}}$) and ($\ori{\mbox{\rm b}}$).

($\ori{\mbox{\rm a}}$) This is a consequence of Lemma~\ref{l:jumps}%
($\ori{\mbox{\rm a}}$), since $\riV(F,T_{(-\infty,t]})=\ori{v}(t)$ and,
by virtue of \eq{e:ttos}, $\riV(F,T_{(-\infty,t)})=\ori{v}(t-0)$.

($\ori{\mbox{\rm b}}$) For every $s\in T_{(t,+\infty)}$, the additivity of $\riV$ implies
  \begin{equation*}
\riV(F,T_{[t,+\infty)})=\ori{v}(s)-\ori{v}(t)+\riV(F,T_{[s,+\infty)}).
  \end{equation*}
Passing to the limit as $T\ni s\to t+0$, we get: by Lemma~\ref{l:jumps}%
($\ori{\mbox{\rm b}}$),
  \begin{equation*}
\ori{v}(s)-\ori{v}(t)\to\ori{v}(t+0)-\ori{v}(t)=\limsup_{T\ni s\to t+0}e(F(t),F(s)),
  \end{equation*}
and by \eq{e:ttoi}, $\riV(F,T_{[s,+\infty)})\to\riV(F,T_{(t,+\infty)})$.
\qed\end{pf}

\subsection{Pointwise Selection Principle}

In the proof of our main result (Theorem~\ref{t:BV}), we will need a compactness
theorem in the topology of pointwise convergence (cf.\ Section~\ref{ss:lsc}) for a
(approximating) sequence of multifunctions $F_n:T\to\mbox{c}(M)$, $n\in\Nb$,
reformulated from \cite{JMAA05} for the metric space $(\mbox{c}(M),d_H)$
under consideration.

Given $F:T\to\mbox{c}(M)$, its \emph{modulus of variation\/} is the nondecreasing
sequence $\{\nu_k(F,T)\}_{k=1}^\infty\subset[0,+\infty]$ defined by
  \begin{equation*}
\nu_k(F,T)\!:=\sup\,\sum_{i=1}^k d_H(F(s_i),F(t_i)),
  \end{equation*}
the supremum being taken over all collections $s_1,\dots,s_k,t_1\dots,t_k$ of $2k$
numbers from $T$ such that $s_1\le t_1\le s_2\le t_2\le\dots\le s_{k-1}\le t_{k-1}%
\le s_k\le t_k$. This notion was introduced in \cite{Chantur} in the context of Fourier
series and extensively applied in \cite[Section~11.3]{GNW} for real valued functions.
The general case of metric space valued functions was considered in \cite{JMAA05},
whence we know that $\nu_k(F,T)\le\nV(F,T)$ for all $k\in\Nb$, and
$\nu_k(F,T)\to\nV(F,T)$ as $k\to\infty$.

The following theorem, extending Helly's Selection Theorem, is a \emph{pointwise
selection principle\/} in terms of the modulus of variation:

\begin{thrm}[{\cite[Theorem~1]{JMAA05}}] \label{t:B}
Suppose $\{F_n\}_{n=1}^\infty$ is a sequence of multifunctions
$F_n:T\to\mbox{\rm c}(M)$ such that\/
{\rm(a)} $\limsup_{n\to\infty}\nu_k(F_n,T)=o(k)$%
\footnote{Equality $\mu_k=o(k)$ means as usual that $\lim_{k\to\infty}\mu_k/k=0$.},
and\/ {\rm(b)}~the closure in $\mbox{\rm c}(M)$ of the set $\{F_n(t):n\in\Nb\}$
is compact for all $t\in T$. Then $\{F_n\}_{n=1}^\infty$ admits a subsequence,
which converges in $\mbox{\rm c}(M)$ pointwise on $T$ to a multifunction
$F:T\to\mbox{\rm c}(M)$ such that $\nu_k(F,T)=o(k)$.
\end{thrm}

\section{Proof of the Main Result} \label{s:pro}

\begin{pfBV}
We prove only item (a), item (b) being proved similarly with corresponding modifications
(`to the left'). We divide the proof into six steps for clarity.
Recall that $a=\inf T\in T$, and we set $b:=\sup T$.
In the first four steps, we prove the theorem in the case when $t_0=a$
and $b\in T$ (so that $T\subset[t_0,b]$ is bounded), in step~5---when $t_0=a$ and
$b\notin T$, and in step~6---when $t_0\!>\!a$ and $T$ is arbitrary.
We employ several ideas from \cite{JMAA,MatSb}, \cite[Sections~5,~10]{Sovae}.

\smallbreak
{\it Step~1}. Suppose $T\subset[t_0,b]$ and $t_0,b\in T$. Since $\riV(F,T)$ is finite
and the variation function of $F$ to the right $\ori{v}:T\to[0,+\infty)$ is nondecreasing
on~$T$, the set of its points of discontinuity on $T$ is at most countable. Denote by
$T_F$ the set of points $t\in T$, which are left limit points of $T$ such that
$\ori{v}(t)=\ori{v}(t-0)$. It follows that $T\setminus T_F$ is at most countable and,
by Lemma~\ref{l:jumps}${\rm(\ori{\mbox{\rm a}})}$,
  \begin{equation} \label{e:FsFt}
\mbox{if $t\in T_F$, then $e(F(s),F(t))\to0$ as $T\ni s\to t-0$.}
  \end{equation}
Furthermore, the set $T_0$ of left isolated points of $T$
(i.e., $t\in T$ such that $T_{(t-\vep,t)}=\es$ for some $\vep>0$) is also at most
countable (in fact, the intervals of ``emptiness from the left'', corresponding to
different left isolated points of $T$, are disjoint, and each such interval contains
a rational number). Let $Q$ denote an at most countable dense subset of $T$. We set
  \begin{equation*}
S\!:=\{t_0,b\}\cup(T\setminus T_F)\cup T_0\cup Q
  \end{equation*}
and note that $S\subset T$ is dense in $T$ and at most countable. With no loss
of generality we may assume that $S$ is countable, say, $S=\{t_i\}_{i=0}^\infty$.
Given $n\in\Nb$, the set $\pi_n=\{t_i\}_{i=0}^{n-1}\cup\{b\}$ is a partition of $T$.
Ordering the points in $\pi_n$ in ascending order and denoting the resulting (ordered)
partition of $T$ by $\pi_n=\{t_i^n\}_{i=0}^n$, we get
  \begin{eqnarray}
&\mbox{$t_0=t_0^n<t_1^n<\dots<t_{n-1}^n<t_n^n=b$, \,and}&\label{e:op}\\[2pt]
&\mbox{$\forall\,t\in S$ $\exists\,n_0=n_0(t)\in\Nb$ such that $t\in\pi_n$
  for all $n\ge n_0$.}& \label{e:noN}
  \end{eqnarray}

{\it Step~2}. Let us construct an approximating sequence for the desired set-valued
selector of $F$. In order to do this, we need the following observation from
\cite[assertion (10.2)]{Sovae}: given $X,Y\in\mbox{c}(M)$,
  \begin{equation} \label{e:proj}
\mbox{there is $Y'\in\mbox{c}(M)$ such that $Y'\subset Y$ and
$d_H(X,Y')\le e(X,Y)$.}
  \end{equation}
In fact, it suffices to define $Y'$ as the \emph{metric projection\/} of $X$ onto~$Y$:
  \begin{equation} \label{e:Pr}
Y'=\mbox{\rm Pr}_Y X:=\{y\in Y:\mbox{there is $x\in X$ such that $d(x,y)=d(x,Y)$}\}.
  \end{equation}

First, given $n\in\Nb$, we define sets $Y_i^n\in\mbox{c}(M)$, $i=0,1,\dots,n$,
inductively as follows. Setting $X=X_0$ and $Y=F(t_0)$ in \eq{e:proj}, we choose
$Y_0:=Y'\in\mbox{c}(M)$ such that $Y_0\subset F(t_0)$ and
$d_H(X_0,Y_0)\le e(X_0,F(t_0))$. We put $Y_0^n:=Y_0$ (for all $n\in\Nb$).
Now, suppose $i\in\{1,\dots,n\}$ and the set $Y_{i-1}^n\in\mbox{c}(M)$ such that
$Y_{i-1}^n\subset F(t_{i-1}^n)$ is already chosen. Then, we put $X=Y_{i-1}^n$
and $Y=F(t_i^n)$ in \eq{e:proj} and pick $Y_i^n:=Y'\in\mbox{c}(M)$ such that
$Y_i^n\subset F(t_i^n)$ and
  \begin{equation} \label{e:Yin}
d_H(Y_{i-1}^n,Y_i^n)\le e(Y_{i-1}^n,F(t_i^n))\le e(F(t_{i-1}^n),F(t_i^n)).
  \end{equation}

The approximating sequence $\Gamma_n:T\to\mbox{c}(M)$, $n\in\Nb$, is defined as
a sequence of set-valued step functions of the form (cf. \cite[equation~(10.3)]{Sovae}):
  \begin{align}
\Gamma_n(t_i^n)&:=Y_i^n\quad\mbox{for all \,$i=0,1,\dots,n;$}\label{e:4+}\\[2pt]
\Gamma_n(t)&:=Y_{i-1}^n\,\,\,\mbox{for all $t\in T\cap(t_{i-1}^n,t_i^n)$
  and $i=1,\dots,n$}\label{e:re}
  \end{align}
(if $T\cap(t_{i-1}^n,t_i^n)=\es$, then $\Gamma_n$ is left undefined
on $(t_{i-1}^n,t_i^n)$). Clearly,
  \begin{equation*}
\Gamma_n(t_0)=\Gamma_n(t_0^n)=Y_0^n=Y_0\subset F(t_0)\quad
\mbox{for all \,$n\in\Nb$.}
  \end{equation*}
Moreover, by the additivity of $\nV$ and \eq{e:Yin}, we have
  \begin{align}
\nV(\Gamma_n,T)&=\sum_{i=1}^n\nV(\Gamma_n,T\cap[t_{i-1}^n,t_i^n])
  =\sum_{i=1}^n d_H(Y_{i-1}^n,Y_i^n)\nonumber\\
&\le\sum_{i=1}^n e(F(t_{i-1}^n),F(t_i^n))\le\riV(F,T)
  \quad\mbox{for all $n\in\Nb$.} \label{e:gn}
  \end{align}

{\it Step~3}. Let us show that $\{\Gamma_n\}_{n=1}^\infty$ satisfies the assumptions
of Theorem~\ref{t:B}. By virtue of \eq{e:gn}, we get
  \begin{equation*}
\limsup_{n\to\infty}\nu_k(\Gamma_n,T)\le\riV(F,T)\quad\mbox{for all \,$k\in\Nb$,}
  \end{equation*}
and so, condition (a) in Theorem~\ref{t:B} is satisfied. Now, we verify condition~(b).

If $t\in S$, then, by \eq{e:noN}, there is $n_0=n_0(t)\in\Nb$ such that $t\in\pi_n$ for
all $n\ge n_0$. So, for each $n\ge n_0$ there is $i=i(t,n)\in\{0,1,\dots,n\}$ such that
$t=t_i^n$. The definition of $\Gamma_n$ implies
  \begin{equation} \label{e:Gat}
\Gamma_n(t)=\Gamma_n(t_i^n)=Y_i^n\subset F(t_i^n)=F(t)\quad
\mbox{for all \,$n\ge n_0$.}
  \end{equation}
In other words, $\{\Gamma_n(t)\}_{n=n_0}^\infty\subset\mbox{c}(F(t))$. Since
$F(t)$ is a compact subset of $M$, it follows from \cite[II.1.4]{CV} that
$\mbox{c}(F(t))$ is a compact subset of $\mbox{c}(M)$, which implies the desired
property for $\{F_n(t)\}_{n=1}^\infty$.

Suppose now that $t\in T\setminus S$. We have $t\in T_F\cap(t_0,b)$, i.e., by \eq{e:FsFt},
  \begin{equation} \label{e:*7}
e(F(s),F(t))\to0\quad\mbox{as}\quad T_{(t_0,b)}\ni s\to t-0,
  \end{equation}
where $t$ is a left limit point of $T$, and so,
there is a sequence $\tau_k\in T$, $\tau_k<t$, $k\in\Nb$, such that $\tau_k\to t$
as $k\to\infty$. Since $S$ is dense in $T$, given $k\in\Nb$, there is $s_k\in S$ such
that $|s_k-\tau_k|<t-\tau_k$, which implies $s_k<t$, and $s_k\to t$ as $k\to\infty$.
From \eq{e:noN}, for each $k\in\Nb$ choose a number $n_0(k)\in\Nb$ (depending
also on~$t$) such that $s_k\in\pi_n$ for all $n\ge n_0(k)$. We may assume (arguing
inductively) that the sequence $\{n_0(k)\}_{k=1}^\infty$ is strictly increasing.
Given $k\in\Nb$ and $n\ge n_0(k)$, since $s_k\in\pi_n$, there is a number
$j(k,n)\in\{0,1,\dots,n-1\}$ such that $s_k=t_{j(k,n)}^n$ and, since $s_k<t$, there is
(unique) $i(k,n)\in\{j(k,n),\dots,n-1\}$ such that
  \begin{equation} \label{e:ijkn}
s_k=t_{j(k,n)}^n\le t_{i(k,n)}^n<t<t_{i(k,n)+1}^n.
  \end{equation}
By \eq{e:4+}, \eq{e:re}, and \eq{e:ijkn}, we find
  \begin{equation} \label{e:Fikn}
\Gamma_n(t)=Y_{i(k,n)}^n\subset F(t_{i(k,n)}^n)\quad\mbox{for all \,$k\in\Nb$
  \,and \,$n\ge n_0(k)$.}
  \end{equation}
Setting $n:=n_0(k)$ and $p_k:=t_{i(k,n_0(k))}^{n_0(k)}$ in \eq{e:Fikn}, we have
$\Gamma_{n_0(k)}(t)\subset F(p_k)$ for all $k\in\Nb$, where, by virtue of \eq{e:ijkn}
and property $s_k\to t$ as $k\to\infty$,
  \begin{equation} \label{e:pkt}
\mbox{$p_k<t$ \,and \,$p_k\to t$ \,as \,$k\to\infty$.}
  \end{equation}
Applying \eq{e:proj}, for each $k\in\Nb$ pick $Y_k(t)\in\mbox{c}(M)$ such that
$Y_k(t)\subset F(t)$ and
  \begin{equation} \label{e:Gok}
d_H(\Gamma_{n_0(k)},Y_k(t))\le e(\Gamma_{n_0(k)},F(t))\le e(F(p_k),F(t)).
  \end{equation}
It follows from \eq{e:*7} and \eq{e:pkt} that the right-hand side in \eq{e:Gok} tends
to zero as $k\to\infty$. By the compactness of $\mbox{c}(F(t))$, we may assume
(passing to a subsequence of $\{Y_k(t)\}_{k=1}^\infty$ if necessary) that
$d_H(Y_k(t),Y(t))\to0$ as $k\to\infty$ for some $Y(t)\in\mbox{c}(M)$ such that
$Y(t)\subset F(t)$. Thus,
  \begin{equation*}
d_H(\Gamma_{n_0(k)}(t),Y(t))\!\le\! d_H(\Gamma_{n_0(k)}(t),Y_k(t))
  \!+\!d_H(Y_k(t),Y(t))\!\to\!0\quad\mbox{as \,$k\!\to\!\infty$,}
  \end{equation*}
and so, (a subsequence of) the subsequence $\{\Gamma_{n_0(k)}\}_{k=1}^\infty$
of $\{\Gamma_n(t)\}_{n=1}^\infty$ converges in $\mbox{c}(M)$. This finishes
the proof of compactness of $\ov{\{\Gamma_n(t)\}_{n=1}^\infty}$ (the closure
being taken in $\mbox{c}(M)$).

\smallbreak
{\it Step~4}. By Theorem~\ref{t:B}, there are a subsequence of
$\{\Gamma_n\}_{n=1}^\infty$, which we denote by
$\{\Gamma_{l(n)}\}_{n=1}^\infty$ with strictly increasing $l:\Nb\to\Nb$, and
a multifunction $\Gamma:T\to\mbox{c}(M)$ such that
$d_H(\Gamma_{l(n)}(t),\Gamma(t))\to0$ as $n\to\infty$ for all $t\in T$.

Let us show that $\Gamma$ is a set-valued selector of $F$ on $T$. It is clear from
\eq{e:4+} and \eq{e:op} that $\Gamma(t_0)=Y_0\subset F(t_0)$, and so (cf. Step~2),
  \begin{equation} \label{e:Xo}
d_H(X_0,\Gamma(t_0))\le e(X_0,F(t_0)).
  \end{equation}

If $t\in S$, \eq{e:noN} implies $t\in\pi_n$ for some $n_0=n_0(t)\in\Nb$ and all $n\ge n_0$.
Since $l(n)\ge l(n_0)\ge n_0$ for $n\ge n_0$, it follows from \eq{e:Gat} that
$\Gamma_{l(n)}(t)\subset F(t)$, and so
(cf.\ properties (a) and (b) of $e$ in Section~\ref{s:mr}),
  \begin{align*}
e(\Gamma(t),F(t))&\le e(\Gamma(t),\Gamma_{l(n)}(t))
  +e(\Gamma_{l(n)}(t),F(t)) \\[2pt]
&=e(\Gamma(t),\Gamma_{l(n)}(t))\le d_H(\Gamma(t),\Gamma_{l(n)}(t))
  \to0\quad\mbox{as \,$n\to\infty$.}
  \end{align*}
This gives $e(\Gamma(t),F(t))=0$ implying $\Gamma(t)\subset F(t)$.

Now suppose that $t\in T\setminus S$. Given $k\in\Nb$, let $n_0(k)$ be the number
(also depending on $t$) from Step~3 such that $s_k\in\pi_n$ for all $n\ge n_0(k)$.
Hence, assertions \eq{e:ijkn} and \eq{e:Fikn} still hold. Setting
$n=N(k):=l(n_0(k))\ge n_0(k)$ and $q_k:=t_{i(k,N(k))}^{N(k)}$ in
\eq{e:ijkn} and \eq{e:Fikn}, we find $\Gamma_{N(k)}\subset F(q_k)$, where
$q_k<t$, and $q_k\to t$ as $k\to\infty$. For each $k\in\Nb$, thanks to \eq{e:proj},
let $Z_k(t)\in\mbox{c}(M)$ be such that $Z_k(t)\subset F(t)$ and
  \begin{equation*}
d_H(\Gamma_{N(k)}(t),Z_k(t))\le e(\Gamma_{N(k)}(t),F(t))\le e(F(q_k),F(t)).
  \end{equation*}
Noting that $e(Z_k(t),F(t))=0$, we get
  \begin{equation*}
e(\Gamma(t),F(t))\le e(\Gamma(t),\Gamma_{N(k)}(t))
  +e(\Gamma_{N(k)}(t),Z_k(t))+e(Z_k(t),F(t))
  \end{equation*}
with the right-hand side tending to zero as $n\!\to\!\infty$. Hence
\mbox{$e(\Gamma(t),F(t))\!=\!0$}, and so, $\Gamma(t)\subset F(t)$.
Thus, we have shown that $\Gamma(t)\subset F(t)$ for all $t\in T$.

The lower semicontinuity of $\nV$ and \eq{e:gn} imply
  \begin{equation*}
\nV(\Gamma,T)\le\liminf_{n\to\infty}\nV(\Gamma_{l(n)},T)\le\riV(F,T).
  \end{equation*}
Since $T_{[a,t_0)}\!=\!\es$ (recall that $t_0\!=\!a$),
$T_{[t_0,+\infty)}\!=\!T_{[t_0,b]}\!=\!T$, and \mbox{$J_a(\Gamma,t_0)\!=\!0$},
this proves Theorem~\ref{t:BV}(a) in the case when $t_0=a$ and $b=\sup T\in T$.

\smallbreak
{\it Remark}. Note that if $X_0\subset F(t_0)$, we have, by \eq{e:Xo},
$\Gamma(t_0)=X_0$.

\smallbreak
{\it Step~5}. Suppose $t_0=a$ and $b=\sup T\notin T$, so that
$T_{[a,t_0)}\!=\!\es$, $T_{[t_0,+\infty)}\!=\!T$, and $J_a(\Gamma,t_0)=0$.
Pick an increasing sequence $\{t_k\}_{k=1}^\infty\subset T$ such that
$t_k\to b$ as $k\to\infty$. Noting that $\riV(F,T_{[t_0,t_1]})\le\riV(F,T)<+\infty$
and applying Steps 1--4 to $F$ on $T_{[t_0,t_1]}$, we get a set-valued selector
$\Gamma_0:T_{[t_0,t_1]}\to\mbox{c}(M)$ of $F$ on $T_{[t_0,t_1]}$ such that
  \begin{equation*}
d_H(X_0,\Gamma_0(t_0))\le e(X_0,F(t_0))\quad\mbox{and}\quad
  \nV(\Gamma_0,T_{[t_0,t_1]})\le\riV(F,T_{[t_0,t_1]}).
  \end{equation*}
Inductively, if $k\in\Nb$, and a set-valued selector $\Gamma_{k-1}$ of $F$ on
$T_{[t_{k-1},t_k]}$ is already chosen, we note that
$\riV(F,T_{[t_k,t_{k+1}]})\le\riV(F,T)<+\infty$, again apply Steps 1--4, and find
a set-valued selector $\Gamma_k$ of $F$ on $T_{[t_k,t_{k+1}]}$ such that
$\Gamma_k(t_k)=\Gamma_{k-1}(t_k)$ and
$\nV(\Gamma_k,T_{[t_k,t_{k+1}]})\le\riV(F,T_{[t_k,t_{k+1}]})$.
Since $t_k\to b$ as $k\to\infty$, and $b\notin T$, we have
$T=\bigcup_{k=0}^\infty T_{[t_k,t_{k+1}]}$, so if $t\in T$ and $t\in T_{[t_k,t_{k+1}]}$
for some $k\in\{0\}\cup\Nb$, we set $\Gamma(t):=\Gamma_k(t)$. Clearly,
$\Gamma:T\to\mbox{c}(M)$ is a well-defined set-valued selector of $F$ on $T$,
inequality \eq{e:Xo} holds, and, by \eq{e:ttos} and the additivity of $\nV$ and~$\riV$,
  \begin{align*}
\nV(\Gamma,T)&=\lim_{n\to\infty}\nV(\Gamma,T_{(-\infty,t_n]})
  =\lim_{n\to\infty}\sum_{k=0}^{n-1}\nV(\Gamma_k,T_{[t_k,t_{k+1}]})\\
&\le\lim_{n\to\infty}\sum_{k=0}^{n-1}\riV(F,T_{[t_k,t_{k+1}]})
  =\lim_{n\to\infty}\riV(F,T_{(-\infty,t_n]})=\riV(F,T).
  \end{align*}

{\it Step~6}. Finally, suppose $t_0>a$. Since $\riV(F,T_{[a,t_0)})$ and
$\riV(F,T_{[t_0,+\infty)})$ do not exceed $\riV(F,T)<+\infty$, and $X_0\in\mbox{c}(M)$,
we apply Steps 1--5 twice: first, to $F$ on $T_{[a,t_0)}$ with arbitrary
$K_0\in\mbox{c}(M)$ in order to obtain a set-valued selector $\Gamma_1$ of $F$
on $T_{[a,t_0)}$ such that
  \begin{equation*}
d_H(K_0,\Gamma_1(a))\le e(K_0,F(a))\quad\mbox{and}\quad
  \nV(\Gamma_1,T_{[a,t_0)})\le\riV(F,T_{[a,t_0)}),
  \end{equation*}
and, second, to $F$ on $T_{[t_0,+\infty)}$ in order to obtain a set-valued selector
$\Gamma_2$ of $F$ on $T_{[t_0,+\infty)}$ such that
  \begin{equation*}
d_H(X_0,\Gamma_2(t_0))\le e(X_0,F(t_0))\quad\mbox{and}\quad
  \nV(\Gamma_2,T_{[t_0,+\infty)})\le\riV(F,T_{[t_0,+\infty)}).
  \end{equation*}
Given $t\in T$, we set $\Gamma(t):=\Gamma_1(t)$ if $t\in T_{[a,t_0)}$, and
$\Gamma(t):=\Gamma_2(t)$ if $t\in T_{[t_0,+\infty)}$.
\linebreak
Clearly, $\Gamma:T\to\mbox{c}(M)$ is a set-valued selector of $F$ on $T$, inequality
\eq{e:Xo} holds, $\nV(\Gamma,T_{[a,t_0)})\le\riV(F,T_{[a,t_0)})$, and
$\nV(\Gamma,T_{[t_0,+\infty)})\le\riV(F,T_{[t_0,+\infty)})$. Furthermore,
since $\riV$ is additive, the second inequality in \eq{e:s1t0} holds (recall that $a=\min T$):
  \begin{align*}
\riV(F,T_{[a,t_0)})+\riV(F,T_{[t_0,+\infty)})&\le
  \riV(F,T_{[a,t_0]})+\riV(F,T_{[t_0,+\infty)})=\\[2pt]
&=\riV(F,T_{[a,+\infty)})=\riV(F,T).
  \end{align*}

Let us prove the first inequality in \eq{e:s1t0}.

Suppose $s_0\in T_{[a,t_0)}$. Hence $s_0=\max T_{[a,t_0)}<t_0$, and so,
$T_{[a,t_0)}=T_{[a,s_0]}$ and $T_{[s_0,t_0]}=\{s_0,t_0\}$ (two-point set).
By the additivity of~$\nV$,
  \begin{align*}
\nV(\Gamma,T_{[a,t_0]})&=\nV(\Gamma,T_{[a,s_0]})+\nV(\Gamma,T_{[s_0,t_0]})
  \\[2pt]
&=\nV(\Gamma_1,T_{[a,t_0)})+d_H(\Gamma_1(s_0),\Gamma_2(t_0)),
  \end{align*}
which implies
  \begin{align*}
\nV(\Gamma,T)&=\nV(\Gamma,T_{[a,t_0]})+\nV(\Gamma,T_{[t_0,+\infty)})\\[2pt]
&=\nV(\Gamma_1,T_{[a,t_0)})+d_H(\Gamma_1(s_0),\Gamma_2(t_0))
  +\nV(\Gamma_2,T_{[t_0,+\infty)})\\[2pt]
&\le\riV(F,T_{[a,t_0)})+d_H(\Gamma(s_0),\Gamma(t_0))+\riV(F,T_{[t_0,+\infty)}).
  \end{align*}
This proves inequality \eq{e:s1t0} with 
$J_a(\Gamma,t_0)=d_H(\Gamma(s_0),\Gamma(t_0))$ from \eq{e:J}.

Now, suppose $s_0\notin T_{[a,t_0)}$. Note that $s_0$ is a left limit point of
$T_{[a,s_0)}$ (and so, of  $T_{[a,t_0)}$ and $T$ as well).
In fact, by the definition of $s_0$, $t<s_0$
for all $t\in T_{[a,s_0)}$, and, given $\vep>0$, there is $t_\vep\in T_{[a,s_0)}$
such that $s_0-\vep<t_\vep$. Hence $s_0-\vep<t_\vep<s_0$, i.e.,
$T_{[a,s_0)}\cap(s_0-\vep,s_0)\ne\es$ (and a fortiori $T\cap(s_0-\vep,s_0)\ne\es$).

Let us show that the limit in the right-hand side of \eq{e:J} exists in $[0,+\infty)$.
Since $\Gamma=\Gamma_1$ on $T_{[a,t_0)}$ and
$\nV(\Gamma_1,T_{[a,t_0)})\le\riV(F,T_{[a,t_0)})$, given $t,t'\in T_{[a,s_0)}$
with $t\le t'$, the additivity of $\nV$ implies (note again that $a=\min T$)
  \begin{align*}
0&\le\nV(\Gamma,T_{[t,t_0]})-\nV(\Gamma,T_{[t',t_0]})=\nV(\Gamma,T_{[t,t']})
  \\[2pt]
&=\nV(\Gamma_1,T_{(-\infty,t']})-\nV(\Gamma_1,T_{(-\infty,t]}).
  \end{align*}
By \eq{e:ttos}, the right-hand side here tends to
$\nV(\Gamma_1,T_{[a,s_0)})-\nV(\Gamma_1,T_{[a,s_0)})=0$ as $T\ni t,t'\to s_0-0$,
and so, Cauchy's criterion yields the existence of the~limit.

Applying the additivity of $\nV$ once again, we get
  \begin{equation} \label{e:star}
\nV(\Gamma,T_{[a,t_0]})=\nV(\Gamma,T_{[a,t]})+\nV(\Gamma,T_{[t,t_0]})
\quad\mbox{for all \,$t\in T_{[a,s_0)}$.}
  \end{equation}
Noting that $s_0\notin T_{[a,t_0)}$ implies $T_{[a,s_0)}=T_{[a,t_0)}$, by virtue
of \eq{e:ttos}, the limit of the first term in \eq{e:star} as $T\ni t\to s_0-0$ is equal to
  \begin{align*}
\lim_{T\ni t\to s_0-0}\nV(\Gamma,T_{[a,t]})&=
  \lim_{T\ni t\to s_0-0}\nV(\Gamma_1,T_{(-\infty,t]})=
  \nV(\Gamma_1,T_{[a,s_0)})\\[2pt]
&=\nV(\Gamma_1,T_{[a,t_0)})\le\riV(F,T_{[a,t_0)}).
  \end{align*}
Taking into account \eq{e:J}, it follows from \eq{e:star} that
  \begin{align*}
\nV(\Gamma,T)&=\nV(\Gamma,T_{[a,t_0]})+\nV(\Gamma,T_{[t_0,+\infty)})\\[2pt]
&=\lim_{T\ni t\to s_0-0}\!\nV(\Gamma,T_{[a,t]})+
  \lim_{T\ni t\to s_0-0}\!\nV(\Gamma,T_{[t,t_0]})+
  \nV(\Gamma_2,T_{[t_0,+\infty)})\\
&=\nV(\Gamma_1,T_{[a,t_0)})+J_a(\Gamma,t_0)+\nV(\Gamma_2,T_{[t_0,+\infty)})
  \\[2pt]
&\le\riV(F,T_{[a,t_0)})+J_a(\Gamma,t_0)+\riV(F,T_{[t_0,+\infty)}).
  \end{align*}
This proves the first inequality in \eq{e:s1t0}.

It remains to show that if $s_0=t_0$, then $J_a(\Gamma,t_0)$ is the left jump of
$\Gamma$ at~$t_0$.\linebreak
Noting that, by the additivity of $\nV$,
  \begin{equation*}
\nV(\Gamma,T_{[t,t_0]})=\ori{v}(t_0)-\ori{v}(t)\quad\mbox{for all \,$t\in T_{[a,t_0)}$,}
  \end{equation*}
and, by \eq{e:VFtm},
  \begin{equation*}
\ori{v}(t_0)=\ori{v}(t_0-0)+\lim_{T\ni t\to t_0-0}d_H(\Gamma(t),\Gamma(t_0)),
  \end{equation*}
and passing to the limit as $T\ni t\to t_0-0$, we get
  \begin{align*}
J_a(\Gamma,t_0)&=\lim_{T\ni t\to t_0-0}\nV(\Gamma,T_{[t,t_0]})
  =\ori{v}(t_0)-\lim_{T\ni t\to t_0-0}\ori{v}(t)\\
&=\ori{v}(t_0-0)+\lim_{T\ni t\to t_0-0}d_H(\Gamma(t),\Gamma(t_0))-\ori{v}(t_0-0)\\
&=\lim_{T\ni t\to t_0-0}d_H(\Gamma(t),\Gamma(t_0)).
  \end{align*}

This completes the proof of Theorem~\ref{t:BV}.
\qed\end{pfBV}

Now, we are in a position to show that Theorem~\ref{t:BV} implies Theorem~\ref{t:A}.

\begin{pfA}
Setting $T_+=T_{[t_0,+\infty)}$ and $T_-=T_{(-\infty,t_0]}$, we have
  \begin{equation*}
\riV(F,T_+)\le\nV(F,T_+)\le\nV(F,T)\quad\mbox{and}\quad
  \leV(F,T_-)\le\nV(F,T_-)\le\nV(F,T).
  \end{equation*}
Noting that $t_0\in T_+\cap T_-$ and applying Theorem~\ref{t:BV}(a) to $F$ on $T_+$
and Theorem~\ref{t:BV}(b) to $F$ on $T_-$, we obtain a set-valued selector
$\Gamma_+:T_+\to\mbox{c}(M)$ of $F$ on $T_+$ and a set-valued selector
$\Gamma_-:T_-\to\mbox{c}(M)$ of $F$ on $T_-$ such that
  \begin{eqnarray*}
&d_H(X_0,\Gamma_+(t_0))\!\le\! e(X_0,F(t_0))\quad\!\mbox{and}\quad\!
  \nV(\Gamma_+,T_+)\!\le\!\riV(F,T_+),&\\[2pt]
&d_H(\Gamma_+(t_0),\Gamma_-(t_0))\!\le\! e(\Gamma_+(t_0),F(t_0))\!=\!0\quad\!
  \mbox{and}\!\quad\nV(\Gamma_-,T_-)\!\le\!\leV(F,T_-).&
  \end{eqnarray*}
Noting that $\Gamma_-(t_0)=\Gamma_+(t_0)$, we set $\Gamma(t):=\Gamma_+(t)$
if $t\in T_+$, and $\Gamma(t):=\Gamma_-(t)$ if $t\in T_-\setminus\{t_0\}$.
Clearly, $\Gamma:T\to\mbox{c}(M)$ is a set-valued selector of $F$ on $T$,
$d_H(X_0,\Gamma(t_0))\le e(X_0,F(t_0))$ and, by the additivity of~$\nV$,
  \begin{align}
\nV(\Gamma,T)&=\nV(\Gamma,T_+)+\nV(\Gamma,T_-)=
  \nV(\Gamma_+,T_+)+\nV(\Gamma_-,T_-)\nonumber\\[1pt]
&\le\riV(F,T_+)+\leV(F,T_-) \label{RR}\\[2pt]
&\le\nV(F,T_+)+\nV(F,T_-)=\nV(F,T). \nonumber
  \end{align}
This finishes the proof of Theorem~\ref{t:A}. Note that we have shown a little
bit more: inequality \eq{RR} holds provided $\riV(F,T_+)$ and $\leV(F,T_-)$ are finite.
\qed\end{pfA}

\begin{rem} \label{r:single}
If $X_0=\{x_0\}\subset M$, a (single-valued) selector of bounded var\-i\-a\-tion
$\Gamma:T\!\to\! M$ of $F$ on $T$ may be chosen such that
\mbox{$d(x_0,\Gamma(t_0))\!\le\! d(x_0,F(t_0))$} and satisfying the rest of assertions
in (a) and (b) of Theorem~\ref{t:BV} (if we replace $d_H$ by $d$ everywhere).
In order to see this, it suffices to pick only one element in the corresponding metric
projection. So (cf. Step~2 in the proof of Theorem~\ref{t:BV}),
choose $y_0\in F(t_0)$ such that $d(x_0,y_0)=d(x_0,F(t_0))$, set $y_0^n:=y_0$,
and if $i\in\{1,\dots,n\}$ and elements $y_{i-1}^n\in F(t_{i-1}^n)$ are already chosen,
pick \mbox{$y_i^n\in F(t_i^n)$} such that $d(y_{i-1}^n,y_i^n)=d(y_{i-1}^n,F(t_i^n))$.
Define $\Gamma_n:T\to M$ (as in \eq{e:4+} and \eq{e:re}) by
$\Gamma_n(t_i^n):=y_i^n$ for $i=0,1,\dots,n$, and $\Gamma_n(t):=y_{i-1}^n$
if $t\in T\cap(t_{i-1}^n,t_i^n)$ and $i=1,\dots,n$. It remains to note
(for $T\subset[t_0,b]$) that
  \begin{align*}
\nV(\Gamma_n,T)&=\nV(\Gamma_n,T\cap[t_{i-1}^n,t_i^n])
  =\sum_{i=1}^n d(y_{i-1}^n,y_i^n)\\
&=\sum_{i=1}^n d(y_{i-1}^n,F(t_i^n))\le\sum_{i=1}^n e(F(t_{i-1}^n),F(t_i^n))
  \le\riV(F,T).
  \end{align*}
In this way, Theorem~\ref{t:BV} above is a generalization of \cite[Theorem~1]{JMAA07},
treating the existence of single-valued selectors on (connected) intervals $T\subset\Rb$.
\end{rem}

\section{Examples} \label{s:exa}

In examples below, we show that all assumptions in Theorem~\ref{t:BV}
are essential.

Let $(\mathbb{B},|\cdot|)$ be a Banach space with norm $|\cdot|$ (e.g.,
$\mathbb{B}=\Rb$), $\mathbb{B}^\Nb$ be the set of all sequences
$x:\Nb\to\mathbb{B}$, and $M=\ell_1(\Nb;\mathbb{B})$ be the (infinite-dimensional)
Banach space of all \emph{summable\/} sequences $x\in\mathbb{B}^\Nb$ equipped
with the norm $\|x\|:=\sum_{i=1}^\infty|x(i)|<+\infty$ and, hence, metric
$d(x,y):=\|x-y\|$ for $x,y\in M$. Fix $u\in\mathbb{B}$ with $|u|=1$
(e.g., $u=1$ in $\Rb$) and, for every $n\in\Nb$, denote by $u_n$ the unit vector
in $M$ defined as usual by $u_n(i)=0$ if $i\ne n$, and $u_n(n)=u$.

\begin{exa} \label{ex1}
A multifunction $F$ on $T\subset\Rb$ with (only) \emph{bounded closed\/} values
in $M$ and $\riV(F,T)<+\infty$ may have no set-valued selectors $\Gamma$ satisfying
  \begin{equation} \label{e:zz}
\nV(\Gamma,T_{[t_0,+\infty)})\le\riV(F,T_{[t_0,+\infty)})\quad\mbox{or}
  \quad d_H(X_0,\Gamma(t_0))\le e(X_0,F(t_0))
  \end{equation}
with $t_0\in T$ and $X_0\in\mbox{c}(M)$. In order to see this, we set
$T:=[0,1]$ and $X:=X_0\cup Y$, where $X_0:=\{u_1\}$ and
$Y:=\{\alpha_nu_n:n\ge2\}$ with $\alpha_n:=1+(1/n)$, and note that $X$ and $Y$
are bounded and closed (but not compact) subsets of $M$, whereas
$X_0\in\mbox{c}(M)$ and $X_0\subset X$.
Define $F:T\to\{X,Y\}\subset\mathcal{P}(M)\setminus\{\es\}$ by
$F(0):=X$ and $F(t):=Y$ if $0<t\le1$. (A similar example was given in
\cite[Example~2]{JMAA} for $F$ with $\nV(F,T)<+\infty$.) We have
$\riV(F,T)=e(X,Y)$, where
  \begin{align*}
e(X,Y)&=\sup_{x\in X}\inf_{y\in Y}d(x,y)=e(X\setminus Y,Y)=e(X_0,Y)\\
&=\inf_{n\ge2}(|u_1|+\alpha_n|u_n|)=1+\inf_{n\ge2}\Bigl(1+\frac1n\Bigr)=2.
  \end{align*}
Since $Y\subset X$, $F$ is nonincreasing on $T$ (Section~\ref{ss:mon}), and so,
$\leV(F,T)=0$.

(a) Suppose $t_0=0$, so that $T_{[t_0,+\infty)}=T$. We have $X_0\subset X=F(0)$,
which implies $e(X_0,F(0))=0$. Let $\Gamma:T\to\mathcal{P}(M)\setminus\{\es\}$
be \emph{any\/} set-valued selector of $F$ on $T$ such that $\Gamma(0)=X_0$.
Since $\es\ne\Gamma(1)\subset F(1)=Y$, $\alpha_nu_n\in\Gamma(1)$ for some
$n\ge2$. It follows that
  \begin{align*}
\nV(\Gamma,T)&\ge d_H(\Gamma(0),\Gamma(1))\ge e(\Gamma(1),\Gamma(0))
  =e(\Gamma(1),X_0)\\[1pt]
&\ge d(\alpha_nu_n,u_1)=\alpha_n+1>2=\riV(F,T).
  \end{align*}

(b) Now, suppose $0<t_0\le1$, and $\Gamma:T\to\mathcal{P}(M)\setminus\{\es\}$
is a set-valued selector of $F$ on $T$. Since $\Gamma(t_0)\subset F(t_0)=Y$,
we have $\alpha_nu_n\in\Gamma(t_0)$ for some $n\ge2$, and so,
  \begin{align*}
d_H(X_0,\Gamma(t_0))&\ge e(\Gamma(t_0),X_0)\ge d(\alpha_nu_n,u_1)
  =\alpha_n+1\\[3pt]
&>2=e(X_0,Y)=e(X_0,F(t_0)).
  \end{align*}

(c) The effect of nonexistence of set-valued selectors in (a) and (b) above is due to
the fact that $\mbox{Pr}_YX_0=\es$ (cf.~\eq{e:Pr}): indeed, if $y\in Y$, then
$y=\alpha_nu_n$ for some $n\ge2$, and so, for $x\in X_0=\{u_1\}$, we have
  \begin{equation*}
d(x,y)=d(u_1,\alpha_nu_n)=1+\alpha_n>2=e(X_0,Y)=d(u_1,Y)=d(x,Y).
  \end{equation*}
\end{exa}

\begin{exa} \label{ex2}
This example is more subtle than Example~\ref{ex1}: even if $F(t)$ is bounded and
closed (but not compact) at a \emph{single\/} point $t\in T$, inequalities \eq{e:zz}
may not hold in Theorem~\ref{t:BV} (this is inspired by \cite[Example~5.2]{Sovae}).

Let $N\in\Nb$, $N\ge2$, be fixed and $\{\alpha_n\}_{n=1}^\infty\subset\Rb$ be
a sequence such that
  \begin{equation} \label{e:C1}
\mbox{$\{|\alpha_n|\}_{n=1}^\infty$ \,is \,strictly \,decreasing}\quad
\mbox{and}\!\!\quad\inf_{n\ge N+1}|\alpha_n|>0
  \end{equation}
(e.g., $\alpha_n=\alpha(n+1)/n$ with $\alpha\ne0$, $n\in\Nb$). We set
$X:=\{\alpha_nu_n:1\le n\le N\}$ and $Y:=\{\alpha_nu_n:n\ge N+1\}$.
Clearly, $X\in\mbox{c}(M)$, while $Y\notin\mbox{c}(M)$ is bounded (by the first
condition in \eq{e:C1}) and closed (by the second condition in \eq{e:C1}) in $M$.
Define $F$ on $T:=[0,1]$ by $F(t)\!:=X$ if $0\le t<1$, and $F(1)\!:=Y$. We have
  \begin{equation*}
\riV(F,T)=e(X,Y)=|\alpha_1|+\inf_{n\ge N+1}|\alpha_n|
  \end{equation*}
and
  \begin{equation*}
\leV(F,T)=e(Y,X)=|\alpha_{N+1}|+|\alpha_N|.
  \end{equation*}

(a) Suppose $0\le t_0<1$ and $X_0:=\{\alpha_1u_1\}$, so that
$T_{[t_0,+\infty)}=[t_0,1]$ and $X_0\subset X=F(t_0)$, which implies
$e(X_0,F(t_0))=0$. Let $\Gamma:T\to\mathcal{P}(M)\setminus\{\es\}$ be a
set-valued selector of $F$ on $T$ such that $\Gamma(t_0)=X_0$. Since
$\Gamma(1)\subset F(1)=Y$, we find $\alpha_nu_n\in\Gamma(1)$ for some
$n\ge N+1$, and so,
  \begin{align*}
\nV(\Gamma,[t_0,1])&\ge d_H(\Gamma(t_0),\Gamma(1))\ge e (\Gamma(1),X_0)
  \ge d(\alpha_nu_n,\alpha_1u_1)=|\alpha_n|+|\alpha_1|\\[2pt]
&>|\alpha_1|+\inf_{i\ge N+1}|\alpha_i|=\riV(F,T)=\riV(F,[t_0,1]).
  \end{align*}

(b) Suppose $t_0=1$ and $X_0:=\{\alpha_1u_1\}$. If
$\Gamma:T\to\mathcal{P}(M)\setminus\{\es\}$ is a set-valued selector of $F$ on $T$,
then $\Gamma(1)\subset F(1)=Y$ implies $\alpha_nu_n\in\Gamma(1)$ for some
$n\ge N+1$, and so,
  \begin{align*}
d_H(X_0,\Gamma(1))&\ge e(\Gamma(1),X_0)\ge d(\alpha_nu_n,\alpha_1u_1)
  =|\alpha_n|+|\alpha_1|\\[2pt]
&>|\alpha_1|+\inf_{i\ge N+1}|\alpha_i|=e(X,Y)\ge e(X_0,Y)=e(X_0,F(t_0)).
  \end{align*}

(c) Clearly, $\mbox{Pr}_XX_0=X_0$. The non-existence of set-valued selectors in (a)
and (b) is again due to the fact that $\mbox{Pr}_YX_0=\es$: in fact, if $y\in Y$,
then $y=\alpha_nu_n$ for some $n\ge N+1$, and so, we have,
for $x\in X_0=\{\alpha_1u_1\}$,
  \begin{equation*}
d(x,y)\!=\!d(\alpha_1u_1,\alpha_nu_n)\!=\!|\alpha_1|+|\alpha_n|
\!>\!\!\inf_{i\ge N+1}(|\alpha_1|+|\alpha_i|)\!=\!\inf_{y\in Y}d(x,y)\!=\!d(x,Y).
  \end{equation*}

(d) We claim that Theorem~\ref{t:BV}(b) holds with $t_0=1$ and $\Gamma(1)=X_0$
(except that $\Gamma(1)\in\mbox{c}(M)$) for every nonempty $X_0\subset Y=F(1)$.

First, observe that $\mbox{Pr}_XY_0=\{\alpha_Nu_N\}$ for every
$\es\ne Y_0\subset Y$ (recall that $X\in\mbox{c}(M)$, and $\mbox{Pr}_XY_0$ is the
set of those $x\in X$, for which $d(y_0,x)=d(y_0,X)$ for some $y_0\in Y_0$).
To see this, we set $n_0:=\min\{n\ge N+1:\alpha_nu_n\in Y_0\}$ for $Y_0\subset Y$.
If $y_0\in Y_0$, we have $y_0=\alpha_nu_n$ for some $n\ge n_0$,
  \begin{equation*}
d(y_0,X)=\inf_{x\in X}d(y_0,x)=\min_{1\le i\le N}(|\alpha_n|+|\alpha_i|)
  =|\alpha_n|+|\alpha_N|,\\[2pt]
  \end{equation*}
and
  \begin{equation*}
\mbox{$d(y_0,x)=|\alpha_n|+|\alpha_i|$ \,if \,$x=\alpha_iu_i\in X$ for some $1\le i\le N$.}
  \end{equation*}
If $n=n_0$ and $i=N$, we find $y_0=\alpha_{n_0}u_{n_0}\in Y_0$,
$x=\alpha_Nu_N\in X$, and
  \begin{equation*}
d(y_0,x)=d(\alpha_{n_0}u_{n_0},\alpha_Nu_N)=|\alpha_{n_0}|+|\alpha_N|=d(y_0,X),
  \end{equation*}
which implies $\alpha_Nu_N\in\mbox{Pr}_XY_0$. Now, if $i<N$, then
$|\alpha_i|>|\alpha_N|$, and so, for every $y_0\in Y_0$, we get
  \begin{equation*}
d(y_0,x)=|\alpha_n|+|\alpha_i|>|\alpha_n|+|\alpha_N|=d(y_0,X).
  \end{equation*}
Thus, $\alpha_iu_i\notin\mbox{Pr}_XY_0$ for all $i=1,\dots,N-1$, and we are through.

Taking the above into account, define a set-valued selector $\Gamma$ of $F$ on $T$
by $\Gamma(t):=\{\alpha_Nu_N\}$ if $0\le t<1$, and $\Gamma(1):=X_0$ with
$X_0=Y_0\subset Y$. It remains to note that $n_0\ge N+1$ implies
  \begin{align*}
\leV(\Gamma,T)&=e(Y_0,\{\alpha_Nu_N\})=\sup_{y_0\in Y_0}d(y_0,\alpha_Nu_N)
  =|\alpha_{n_0}|+|\alpha_N|\\
&\le|\alpha_{N+1}|+|\alpha_N|=\leV(F,T).
  \end{align*}
\end{exa}

\begin{exa} \label{ex3}
Making use of an idea from \cite[Example~3.1]{JMAA07}, here we present an example
of a multifunction $F:T\to\mbox{c}(M)$ with $T:=[0,1]$ such that
$\riV(F,T)<+\infty$ and $\leV(F,T)=+\infty$ (thus, Theorem~\ref{t:BV}(a) is
applicable to $F$, whereas Theorem~\ref{t:A} is not).

Let $N\in\Nb$ and $\{\alpha_n\}_{n=1}^\infty\subset(0,+\infty)$ be a
decreasing sequence such that
  \begin{equation} \label{e:Alph}
\lim_{n\to\infty}\alpha_n=0\quad\mbox{and}\quad
\sum_{n=1}^\infty\alpha_{nN}=+\infty
  \end{equation}
(e.g., $\alpha_n=1/n$). Given $n\in\Nb$, we set
  \begin{equation*}
X_n:=\{0\}\cup\{\alpha_iu_i:1\le i\le nN\},\quad\mbox{and}\quad
X_\infty:=\{0\}\cup\{\alpha_iu_i:i\in\Nb\}.
  \end{equation*}
Clearly, $X_n\in\mbox{c}(M)$ for all $n\in\Nb$ and, by the first assumption in \eq{e:Alph},
the set $X_\infty$ is compact as well. Let $\{\tau_n\}_{n=0}^\infty\subset[0,1)$ be a
strictly increasing sequence such that $\tau_0=0$ and $\lim_{n\to\infty}\tau_n=1$.
Define $F:T\to\mbox{c}(M)$ by the rule:
  \begin{equation*}
\mbox{$F(t):=X_n$ \,if\, $\tau_{n-1}\le t<\tau_n$ for all $n\in\Nb$, \,and\,
$F(1):=X_\infty$.}
  \end{equation*}
Since $X_n\subset X_{n+1}\subset X_\infty$ for all $n\in\Nb$, $F$ is nondecreasing
on $T$ (Section~\ref{ss:mon}), and so, $\riV(F,T)=0$. In order to see that
$\leV(F,T)=+\infty$, given $m\in\Nb$, $m\ge2$, and partition
$\pi_m=\{\tau_n\}_{n=0}^{m-1}\cup\{1\}$ of $T=[0,1]$, we find
  \begin{align}
\leV(F,T)&\ge\sum_{n=1}^{m-1}e(F(\tau_n),F(\tau_{n-1}))+e(F(1),F(\tau_{m-1}))
  \nonumber\\
&=\sum_{n=1}^{m-1}e(X_{n+1},X_n)+e(X_\infty,X_m), \label{e:rhs}
  \end{align}
where
  \begin{equation*}
e(X_{n+1},X_n)=\sup_{nN+1\le k\le(n+1)N}\Bigl(|\alpha_k|+
  \inf_{1\le i\le nN}|\alpha_i|\Bigr)=\alpha_{nN+1}+\alpha_{nN}
  \end{equation*}
and
  \begin{equation*}
e(X_\infty,X_m)=\sup_{k\ge mN+1}\Bigl(|\alpha_k|+
  \inf_{1\le i\le mN}|\alpha_i|\Bigr)=\alpha_{mN+1}+\alpha_{mN}.
  \end{equation*}
It follows that the quantity \eq{e:rhs} is equal to
  \begin{equation*}
\sum_{n=1}^m\alpha_{nN+1}+\sum_{n=1}^m\alpha_{nN}\to+\infty\quad
\mbox{as}\quad m\to\infty.
  \end{equation*}
\end{exa}

\begin{exa} \label{ex4}
(a) By Theorem~\ref{t:BV}(a), given $X_0\in\mbox{c}(M)$ such that
$X_0\subset X_1=F(0)$, multifunction $F$ from Example~\ref{ex3} has a constant
set-valued selector $\Gamma:T\to\mbox{c}(M)$ satisfying $\Gamma(0)=X_0$
and $\nV(\Gamma,T)\le\riV(F,T)=0$. However, if $0<t_0\le1$ and $X_0\subset F(t_0)$,
there may be no set-valued selector $\Gamma$ of $F$ on $T$ such that
$\Gamma(t_0)=X_0$ and $\nV(\Gamma,T)\le\riV(F,T)$. In fact, let $t_0=\tau_n$
for some $n\in\Nb$ (cf. Example~\ref{ex3}), so that $F(t_0)=F(\tau_n)=X_{n+1}$.
Suppose now that $X_0\subset X_{n+1}\setminus X_1\subset F(t_0)$,
$\Gamma(t)\subset F(t)$ for all $t\in T$, and $\Gamma(t_0)=X_0$. Since
$\Gamma(0)\subset F(0)=X_1$, we have $\alpha_iu_i\in\Gamma(0)$ for some
$1\le i\le N$, or $0\in\Gamma(0)$ (i.e., possibly, $\alpha_i=0$), and so,
  \begin{align}
\nV(\Gamma,T)&\ge d_H(\Gamma(t_0),\Gamma(0))=d_H(X_0,\Gamma(0))
  \ge e(\Gamma(0),X_0)\nonumber\\[3pt]
&\ge d(\alpha_iu_i,X_0)\ge d(\alpha_iu_i,X_{n+1}\setminus X_1)\nonumber\\[2pt]
&=\min_{N+1\le k\le (n+1)N}(\alpha_i+\alpha_k)\ge\alpha_{(n+1)N}
  >0=\riV(F,T).\label{e:519}
  \end{align}

(b) If $t_0>a=\inf T$ and $\Gamma(t_0)=X_0$, inequality
$\nV(\Gamma,T_{[a,t_0)})\le\riV(F,T_{[a,t_0)})$ in Theorem~\ref{t:BV}(a) cannot
in general be replaced by $\nV(\Gamma,T_{[a,t_0]})\le\riV(F,T_{[a,t_0]})$.
This can be seen from Example~\ref{ex4} and \eq{e:519}:
  \begin{equation*}
\nV(\Gamma,[0,t_0])\ge d_H(\Gamma(t_0),\Gamma(0))\ge\dots\ge
\alpha_{(n+1)N}>0=\riV(F,[0,t_0]).
  \end{equation*}
This observation also makes it explicit that the ``jump'' $J_a(\Gamma,t_0)$ is
essential  in the left-hand side of \eq{e:s1t0}.
\end{exa}

\begin{exa} \label{ex:5.5}
This example is designed for Remark~\ref{p:7}. Let $T:=[1,+\infty)$ and
$F:T\to\mbox{c}(M)$ be given by $F(t)\!:=\!X_n$ if $n\in\Nb$ and $n\le t<n+1$,
where $X_n\!:=\!\{u_i:1\le i\le n\}$. We have $X_n\in\mbox{c}(M)$, and
$F(s)\subset F(t)$ for all $1\le s\le t<+\infty$, and so, $\riV(F,T)=0$. The image
\mbox{$F(T)=\bigcup_{n=1}^\infty X_n=\{u_i:i\in\Nb\}$} is bounded in $M$,
but \emph{not\/} totally bounded (i.e., cannot be covered by a finite number of balls
of arbitrarily small radius). Note that $\leV(F,T)=+\infty$ (consider a partition
$1<2<\dots<m-1<m$ of $T$ with arbitrary $m\in\Nb$ and observe that
$e(X_{n+1},X_n)=2$ for all $n\in\Nb$). This example is easily adapted to the case
when $F$ maps $[a,b)$ or $[a,b]$ into~$\mbox{c}(M)$.
\end{exa}

\section{Functional Inclusion $X(t)\subset F(t,X(t))$} \label{s:Func_Em}

Assuming some interplay of the (uniform) boundedness of directional variations and
(uniform) contractions, we have the following parametrized version of Banach's
Contraction Theorem, extending Theorem~11.4 from~\cite{Sovae}.

\begin{thm} \label{t:FuEm}
Suppose a multifunction $F:T\!\times\!\mbox{\rm c}(M)\!\to\!\mbox{\rm c}(M)$
is such that\/{\rm:}
  \begin{itemize}
\item[{\rm(a)}] there is a nondecreasing bounded function $\varphi:T\to\Rb$ such that
  \begin{equation*}
\mbox{$\!\!\!e(F(s,X),F(t,X))\le\varphi(t)-\varphi(s)$ for all $s,t\in T$, $s\le T$, and
  $X\in\mbox{\rm c}(M);$}
  \end{equation*}
\item[{\rm(b)}] there is a number $0\le\mu<1$ such that
  \begin{equation*}
\mbox{$e(F(t,X),F(t,Y))\le\mu d_H(X,Y)$ for all $t\in T$ and $X,Y\in\mbox{\rm c}(M);$}
  \end{equation*}
\item[{\rm(c)}] there is a multifunction $K:T\to\mbox{\rm c}(M)$ such that
  \begin{equation*}
\mbox{$F(t,X)\subset K(t)$ for all $t\in T$ and $X\in\mbox{\rm c}(M)$.}
  \end{equation*}
  \end{itemize}

If $t_0:=\inf T\in T$ and $X_0\in\mbox{\rm c}(M)$, then there is
$X:T\to\mbox{\rm c}(M)$ such that\\[3pt]
{\rm(i)}~$\nV(X,T)\le\nV(\varphi,T)/(1-\mu)<+\infty;$
{\rm(ii)}~$X(t)\subset F(t,X(t))$ for all $t\in T;$\\[3pt]
{\rm(iii)} $d_H(X_0,X(t_0))\le e\bigl(X_0,F(t_0,X(t_0))\bigr)$.
\smallbreak
In addition, if $X_0\subset F(t_0,X_0)$, then {\rm(iii)} can be replaced by $X(t_0)=X_0$.
\end{thm}

\begin{pf}
First, observe that assumptions (a) and (b) and the triangle inequality for $e$ imply,
for all $s,t\in T$, $s\le t$, and $X,Y\in\mbox{c}(M)$,
  \begin{equation} \label{e:XY}
e(F(s,X),F(t,Y))\le\varphi(t)-\varphi(s)+\mu d_H(X,Y).
  \end{equation}

We set $X_0(t)\!:=\!X_0$ and $F_0(t)\!:=\!F(t,X_0)$ for $t\in T$. We have
\mbox{$F_0:T\to\mbox{c}(M)$}, and assumption (a) and Lemma~\ref{l:char} imply
$\riV(F_0,T)\le\nV(\varphi,T)<+\infty$. By Theorem~\ref{t:BV}(a), there is
$X_1\equiv\Gamma:T=T_{[t_0,+\infty)}\to\mbox{c}(M)$ such that
$X_1(t)\subset F_0(t)$ for all $t\in T$, $d_H(X_0,X_1(t_0))\le e(X_0,F_0(t_0))$, and
$\nV(X_1,T)\le\riV(F_0,T)\le\nV(\varphi,T)$. In what follows we apply the standard
iteration procedure. Setting $F_1(t):=F(t,X_1(t))$ for $t\in T$, we find
$F_1:T\to\mbox{c}(M)$ and, by \eq{e:XY},
  \begin{equation*}
e(F_1(s),F_1(t))\le\varphi(t)-\varphi(s)+\mu d_H(X_1(s),X_1(t))\quad
  \forall\,s,t\in T,\,\,s\le t.
  \end{equation*}
Arguing with partitions of $T$, Lemma~\ref{l:char} implies
  \begin{equation*}
\riV(F_1,T)\le\nV(\varphi,T)+\mu\nV(X_1,T)\le(1+\mu)\nV(\varphi,T).
  \end{equation*}
Applying Theorem~\ref{t:BV} again, we obtain $X_2:T\to\mbox{c}(M)$ such that
$X_2(t)\subset F_1(t)$ for all $t\in T$, $d_H(X_0,X_2(t_0))\le e(X_0,F_1(t_0))$, and
  \begin{equation*}
\nV(X_2,T)\le\riV(F_1,T)\le(1+\mu)\nV(\varphi,T).
  \end{equation*}
If $F_2(t)\!:=\!F(t,X_2(t))$, $t\in T$, then $F_2:T\to\mbox{c}(M)$ and, by \eq{e:XY},
  \begin{equation*}
e(F_2(s),F_2(t))\le\varphi(t)-\varphi(s)+\mu d_H(X_2(s),X_2(t))\quad
  \forall\,s,t\in T,\,\,s\le t,
  \end{equation*}
and so,
  \begin{equation*}
\riV(F_2,T)\le\nV(\varphi,T)+\mu\nV(X_2,T)\le(1+\mu+\mu^2)\nV(\varphi,T).
  \end{equation*}

Arguing by induction, we obtain the sequence $\{X_n\}_{n=1}^\infty$ of multifunctions
$X_n:T\to\mbox{c}(M)$ such that, given $n\in\Nb$,
  \begin{eqnarray}
&X_n(t)\subset F_{n-1}(t)\!:=\!F(t,X_{n-1}(t))\subset K(t)\quad\mbox{for all $t\in T$,}&
  \label{e:Xnt}\\[3pt]
&d_H(X_0,X_n(t_0))\le e(X_0,F_{n-1}(t_0))\!=\!e\bigl(X_0,F(t_0,X_{n-1}(t_0))\bigr),
  \,\,\mbox{and}& \label{e:Hot}\\
&\displaystyle\nV(X_n,T)\le\biggl(\sum_{i=0}^{n-1}\mu^i\biggr)\nV(\varphi,T)
  \le\frac1{1-\mu}\,\nV(\varphi,T).& \label{e:Vft}
  \end{eqnarray}
By \eq{e:Vft}, the sequence $\{X_n\}_{n=1}^\infty$ is of uniformly bounded Jordan
variation with respect to $d_H$, and so, condition (a) in Theorem~\ref{t:B} is satisfied,
and by \eq{e:Xnt}, the closure $\ov{\{X_n(t):n\in\Nb\}}$ in $\mbox{c}(M)$ is compact
for every $t\in T$, and so, condition (b) in Theorem~\ref{t:B} is fulfiled.
By Theorem~\ref{t:B}, a subsequence of $\{X_n\}_{n=1}^\infty$, again denoted by
$\{X_n\}_{n=1}^\infty$, converges in $\mbox{c}(M)$ pointwise on $T$ to a
multifunction $X:T\to\mbox{c}(M)$, i.e., $d_H(X_n(t),X(t))\to0$ as $n\to\infty$
for all $t\in T$.

We are going to verify that $X$ satisfies assertions (i), (ii), and (iii). Assertion (i)
is a consequence of \eq{e:Vft} and the lower semicontinuity \eq{e:Vsem} of~$\nV$.
In order to see that (ii) holds, we make use of the following inequality
(cf. \cite[inequality~(11.7)]{Sovae}), which is valid for all $X,X',Y,Y'\in\mbox{c}(M)$:
  \begin{equation} \label{e:Sog}
|e(X,Y)-e(X',Y')|\le d_H(X,X')+d_H(Y,Y').
  \end{equation}
In fact, given $t\in T$, \eq{e:Xnt} implies $e\bigl(X_n(t),F(t,X_{n-1}(t))\bigr)=0$,
and so, taking into account \eq{e:Sog} and \eq{e:XY}, we get
  \begin{align*}
e(X(t),F(t,X(t)))&\!=\!\bigl|e(X(t),F(t,X(t)))\!-\!e(X_n(t),F(t,X_{n-1}(t)))\bigr|\\[2pt]
&\!\le\! d_H(X(t),X_n(t))\!+\!d_H(F(t,X(t)),F(t,X_{n-1}(t)))\\[2pt]
&\!\le\! d_H(X(t),X_n(t))\!+\!\mu d_H(X(t),X_{n-1}(t))\to0
  \,\,\,\,\mbox{as}\,\,\,n\to\infty.
  \end{align*}
Thus, $e(X(t),F(t,X(t)))=0$, which implies (by properties of~$e$) assertion~(ii).

To establish (iii), we note that (cf. \eq{e:Hot} and (iii))
  \begin{equation*}
|d_H(X_0,X_n(t_0))-d_H(X_0,X(t_0))|\le d_H(X_n(t_0),X(t_0))\to0
  \,\,\,\,\mbox{as}\,\,\,n\to\infty,
  \end{equation*}
and, by virtue of \eq{e:Sog} and assumption (b),
  \begin{align*}
&\qquad\bigl|e(X_0,F(t_0,X_{n-1}(t_0)))-e(X_0,F(t_0,X(t_0)))\bigr|\\[2pt]
\le&\,d_H(F(t_0,X_{n-1}(t_0)),F(t_0,X(t_0)))\le\mu d_H(X_{n-1}(t_0),X(t_0))
  \to0,\,\,n\to\infty.
  \end{align*}
Passing to the limit as $n\to\infty$ in \eq{e:Hot}, we arrive at~(iii).

Finally, suppose $X_0\!\subset\! F(t_0,X_0)$. Hence $X_0\!\subset\! F_0(t_0)$ and
$e(X_0,F_0(t_0))\!=\!0$. From the above,
$d_H(X_0,X_1(t_0))\!\le\! e(X_0,F_0(t_0))\!=\!0$,
and so, $X_1(t_0)\!=\!X_0$. Since
  \begin{equation*}
X_0\subset F_0(t_0)=F(t_0,X_0)=F(t_0,X_1(t_0))=F_1(t_0),
  \end{equation*}
we find from $d_H(X_0,X_2(t_0))\le e(X_0,F_1(t_0))=0$ that $X_2(t_0)=X_0$.
By induction, we deduce from \eq{e:Hot} that $X_n(t_0)=X_0$ for all $n\in\Nb$.
Passing to the lmit as $n\to\infty$ in \eq{e:Hot}, we get
$d_H(X_0,X(t_0))\le e(X_0,F(t_0,X_0))=0$, which yields $X(t_0)=X_0$.
\qed\end{pf}

\begin{rem}
If $F(t,X)=F(t)$ is independent of $X\in\mbox{c}(M)$ (or $\mu=0$), Theorem~%
\ref{t:FuEm} gives back Theorem~\ref{t:BV}(a): we may set $\varphi=\ori{v}_{F}$
and~$K=F$. On the other hand, if $F(t,X)=F(X)$ is independent of $t\in T$ (or
$\varphi\equiv0$), Theorem~\ref{t:FuEm} is a consequence of Banach's Contraction
Theorem (in fact, $F:\mbox{c}(K)\to\mbox{c}(K)$ is a contraction on compact,
hence complete, metric space $(\mbox{c}(K),d_H)$ with $K=K(t_0)$).
\end{rem}

\begin{exa} \label{ex:6.1}
The purpose of this example is to show that assumptions of Theorem~\ref{t:FuEm}
can be fulfiled. Let $M=\mathbb{B}$ be a Banach space with norm $|\cdot|$ and
metric $d(x,y)=|x-y|$, $x,y\in M$, and $K\in\mbox{c}(M)$. Suppose
$\varphi_0:T\to[0,+\infty)$ is nondecreasing and $\mu\!:=\!\sup_{t\in T}\varphi_0(t)<1$.
Define $F:T\times\mbox{c}(M)\to\mbox{c}(M)$ by $F(t,X):=\varphi_0(t)X$ for
$t\in T$ and $X\in\mbox{c}(M)$. We have
  \begin{equation*}
e(F(s,X),F(t,X))\le\bigl(\varphi_0(t)-\varphi_0(s)\bigr)\max_{x\in X}|x|\quad
\mbox{for all $s,t\in T$, $s\le t$,} 
  \end{equation*}
and so, condition (a) in Theorem~\ref{t:FuEm} is satisfied with
$\varphi(t)\!:=\!\varphi_0(t)\max_{x\in K}|x|$, $t\in T$, for all $X\in\mbox{c}(K)$.
Furthermore, given $t\in T$ and $X,Y\in\mbox{c}(K)\subset\mbox{c}(M)$,
  \begin{equation*}
e(F(t,X),F(t,Y))=\varphi_0(t)e(X,Y)\le\mu d_H(X,Y),
  \end{equation*}
and so, condition (b) in Theorem~\ref{t:FuEm} is satisfied. Finally, setting
$K(t):=\varphi_0(t)K$, we find $F(t,X)\subset K(t)$ for all $t\in T$
and $X\in\mbox{c}(K)$.
\end{exa}

\begin{exa}
In Example~\ref{ex:6.1}, we set $\mathbb{B}:=\Rb$, $K:=[0,1]$, and
$\varphi_0(t):=t$ for $t\in T:=[0,1/2]$ (hence $\mu=1/2$, $\varphi=\varphi_0$,
and $\nV(\varphi,T)=1/2$). Define $F$ by
  \begin{equation*}
F(t,X):=(tX)\cup(1-t+tX),\quad t\in T,\,\,X\in\mbox{c}(K).
  \end{equation*}
For instance, if $X=[0,1]$, we have $F(0,X)=\{0,1\}$, $F(1/2,X)=[0,1]$, and if
$0<t<1/2$, then $F(t,X)=[0,t]\cup[1-t,1]$, $F(t,[0,t])=[0,t^2]\cup[1-t,1-t+t^2]$,
$F(t,[1-t,1])=[t-t^2,t]\cup[1-t^2,1]$, and so on. The iterative construction of the
classical Cantor (ternary) set corresponds to $t=1/3$ (e.g., \cite[p.~20]{GK}).

By Theorem~\ref{t:FuEm}, there is $X:[0,1/2]\to\mbox{c}([0,1])$ such that
$\nV(X,[0,1/2])\!\le\!1$, $X(t)\subset(tX(t))\cup(1-t+tX(t))$ for all $t\in[0,1/2]$, and
$X(0)=\{0,1\}$. For every $t\in(0,1/2)$, the compact set $X(t)\subset[0,1]$ is a
Cantor-type perfect set.
\end{exa}

{\it Acknowledgments}.
The article was prepared within the framework of the Academic Fund Program at the
National Research University Higher School of Economics (HSE) in 2017--2018
(grant~no.\,17-01-0050) and by the Russian Academic Excellence Project
\verb|"|5--100\verb|"|.



\end{document}